\documentclass[11pt,a4paper]{article}
\usepackage[utf8]{inputenc}
\usepackage{amsmath}
\usepackage{amsfonts}
\usepackage{amssymb}
\usepackage{bm}
\usepackage{pdflscape}
\usepackage[margin=1in,footskip=0.25in]{geometry}
\usepackage{graphicx}
\usepackage{multirow}
\usepackage{graphicx}
\usepackage{wrapfig}
\usepackage{setspace}
\usepackage{lscape}
\usepackage[table,xcdraw]{xcolor}
\usepackage{float}

\usepackage{longtable}
\usepackage{courier}
\usepackage[printwatermark]{xwatermark}
\begin{document}

\noindent\textbf{\Large{Bias corrected estimators for proportion of true null hypotheses under exponential model: Application of adaptive FDR-controlling in segmented failure data}}\\
\begin{center}
Aniket Biswas\\
Department of Statistics\\
Dibrugarh University\\
Dibrugarh, Assam, India-786004\\
Email:\textit{biswasaniket44@gmail.com}
\end{center}

\begin{center}
Gaurangadeb Chattopadhyay\\
Department of Statistics\\
University of Calcutta\\
Kolkata, West Bengal, India-700019\\
Email:\textit{gcdhstat@gmail.com}
\end{center}

\begin{center}
Aditya Chatterjee\\
Department of Statistics\\
University of Calcutta\\
Kolkata, West Bengal, India-700019\\
Email:\textit{acbustat1@gmail.com}
\end{center}


\section*{Abstract}
Two recently introduced model based bias corrected estimators for proportion of true null hypotheses ($\pi_0$) under multiple hypotheses testing scenario have been restructured for exponentially distributed random observations available for each of the common hypotheses. Based on stochastic ordering, a new motivation behind formulation of some related estimators for $\pi_0$ is given. The reduction of bias for the model based estimators are theoretically justified and algorithms for computing the estimators are also presented. The estimators are also used to formulate a popular adaptive multiple testing procedure. Extensive numerical study supports superiority of the bias corrected estimators. We also point out the adverse effect of using the model based bias correction method without proper assessment of the underlying distribution. A case-study is done with a synthetic dataset in connection with reliability and warranty studies to demonstrate the applicability of the procedure, under a non-Gaussian set up. The results obtained are in line with the intuition and experience of the  subject expert. An intriguing discussion has been attempted to conclude the     article that also indicates the future scope of study.\\

\noindent \textbf{Keywords}: Multiple hypotheses testing, Adaptive Benjamini-Hochberg algorithm, Mean mileage to failure, $p$-value.\\

\noindent \textbf{MS 2010 classification}: 62F99, 62P30, 62N99. \\

\section{Introduction}
Let us consider an empirical Bayesian set-up given in Storey (2002), where $m$ similar but independent hypotheses are to be tested, viz. $H_1, H_2, ..., H_m$. For $H_i=1$, the $i$-th null hypothesis is true and for $H_i=0$, false for any $i\in \{1, 2, ..., m\}$. Thus, $H_i$'s are Bernoulli random variables with success probability $\pi_0\in (0,1)$, the proportion of true null hypotheses. Let $m_0$ be the number of true null hypotheses. Thus, $m_0=\sum_{i=1}^m H_i$ is a binomial random variable with index $m$ and parameter $\pi_0$. Clearly, $H_i$'s and hence $m_0$ remain latent and can never be realized in a given multiple testing scenario. As in case of single hypothesis testing problem, the test statistics $T_1, T_2, ..., T_m$ respectively for $H_1, H_2,..., H_m$ may be observed. For $F_0$ being the common distribution of $T_i|H_i=1$ and $F_1$ being the same for $T_i|H_i=0$, a two-component mixture model for $T_i$ is
\begin{equation}\label{eq:Tmix}
T_i\sim \pi_0\, F_0\,+\,(1-\pi_0)\,F_1\quad\textrm{for all}\quad i=1, 2, ..., m.
\end{equation}
Thus, $\pi_0$ may be thought of as the mixing proportion of the null test statistics with the non-null test statistics when multiple tests are performed. In existing literature $p$-values are considered as test statistics since its use ensures similar nature of critical region, irrespective of the nature of hypotheses framed. Usually, a little abuse of notation is made while denoting $p$-value by $p$ irrespective of whether it is a random variable or a realization on that. The distinction of usage ought to be understood as the situation demands. The marginal density function of $p$-value (Langaas et al. 2005) is
\begin{equation}\label{eq:pmix}
f(p)=\pi_0\,f_0(p)\,+\,(1-\pi_0)\,f_1(p)\quad\textrm{for}\quad0<p<1
\end{equation}
where, $f_0$ and $f_1$ are two p-value densities, respectively under the null and alternative hypotheses. When the tested null is simple and the corresponding test statistic is absolutely continuous, $f_0(p)$ is simply $1$, the density function of a uniform random variable over $(0,1)$ and the $p$-value under the alternative hypothesis is stochastically smaller than the uniform variate. In addition, the density estimation based approaches for estimating $\pi_0$ impose certain restrictions on $f_1$ (Langaas et al. 2005; Guan et al. 2008; Ostrovnaya and Nicolae 2012). Storey's estimator (Storey 2002) is constructed on the basis of a tuning parameter $\lambda\in (0,1)$ such that, $f_1(p)=0$ for $p>\lambda$. This assumption introduces a conservative bias in the estimator that can be corrected or in practice, can be reduced as have been discussed in Cheng et al. (2015). The set-up given therein for the applicability of the Gaussian model based bias correction is discussed in section 2. Biswas (2019) has recently proposed an alternative model based bias corrected estimator for $\pi_0$ under the same set-up. A comparative performance study of both the estimators with simulated microarray datasets has also been provided. Microarray datasets attract the limelight to demonstrate the application of multiple testing procedures owing to their obvious high dimensional nature and associated decision-making questions.\\

The current work deals with a segmented failure dataset, where failure time or some similar entity of a particular component is available for a number of units but the units are operated or tested in different conditions, that may vary over space and time. Thus, the dataset is divided into several segments and the observations are available for each segment. The number of observations per segment (in order of tens or hundreds) might be much less than the number of segments (in order of hundreds or thousands), as the segmentation is done on the basis of time and space among other things. Thus the situation is quite similar to that of microarray datasets where thousands of genes are tested to identify differentially expressed genes based on gene expression levels of two small groups of subjects, viz. treatment group and control group. For segmented failure dataset, similar kind of questions may be raised regarding identification of segment(s) for which the failure patterns of that particular component performs strikingly different (much worse or better) from the average. \\

To answer this question, appropriate hypotheses for each segment are framed and tested simultaneously. While testing a large number of hypotheses, control over the False Discovery Rate (FDR) is desirable and the classical Benjamini-Hochberg algorithm has the ability to do so (Benjamini and Hochberg 1995). A reliable estimate of $\pi_0$ may be used for eliminating the conservativeness present in such Benjamini-Hochberg approach (Benjamini and Yekutieli 2000). Empirical Bayesian interpretation of FDR and controlling the same by estimating it for fixed rejection region requires an estimate of $\pi_0$ (Storey 2002). Craiu and Sun (2008) justifies the equivalence of Storey's q-value approach with adaptive Benjamini-Hochberg algorithm. The authors emphasized that both the adaptive procedure and Storey's approach require a good approximation of $\pi_0$.\\

Gaussian model assumption for failure data is inappropriate, especially when sample size corresponding to each segment or equivalently each test is small and exponential distribution is a reasonable primary model choice. Under this set-up, we modify both the estimators proposed in Cheng et al. (2015) and Biswas (2019) and find that these model based estimators are more efficient than the existing $\pi_0$-estimators in practice. Application of  adaptive Benjamini-Hochberg procedure has the ability  to list the significantly different segments with respect to such time to event or equivalent entity of a certain component in our case study. \\

The remaining part of the article is structured as follows. In section 2 we reproduce Storey's estimator and the recently introduced bias corrected estimators from stochastic ordering approach which ties them in a yarn and may inspire further works in similar line. The next section is devoted to different testing scenarios and useful properties of respective non-null $p$-values. In section 4 we briefly revisit the estimation algorithms and discuss adaptation of the $\pi_0$ estimates to Benjamini-Hochberg algorithm. Section 5 deals with performance comparison of the new estimators with existing ones through extensive simulation experiment. In section 6, a real-life synthetic segmented failure dataset is presented, has been validated for applicability of the proposed methods, analysed to demonstrate the superior performance of adaptive algorithm with the new estimators along with proper justification of the findings. We conclude the article with a  mention of a few limitations of the present work and a glimpse of the future direction of the study.
   
\section{Methods of estimation}
Let $p$ denote a $p$-value corresponding to a simple null hypothesis testing problem with continuous test statistic. Thus, $p$ has the support $(0,1)$. Consider another random variable $V$ on the same support $(0,1)$ with the  distribution function $G$. Then,
\begin{equation}\label{eq:pgreaterV}
P(p\geq V)=\int_{0}^1 f(p)\, G(p)\, dp. 
\end{equation}     
In the following subsections we take different choices for $G$ and motivate different estimators for $\pi_0$ as mentioned in section 1.

\subsection{Storey's bootstrap estimator and related approaches}
Consider $V$ to be degenerate at some $\lambda\in (0,1)$. Thus,
\begin{equation}\label{G1}
G(v)=\begin{cases} 1\quad\textrm{for}\quad v\geq \lambda\\
 0\quad\textrm{for}\quad v< \lambda.\\
 \end{cases}
\end{equation}
Putting (\ref{eq:pmix}) and (\ref{G1}) in (\ref{eq:pgreaterV}), we obtain
\begin{equation}\label{eq:Fbarl}
\bar{F}(\lambda)= \pi_0\,(1-\lambda)\,+\,(1-\pi_0)\, Q(\lambda)
\end{equation}
 where, $F$ is the distribution function of $p$, $\bar{F}=1-F$ and $Q$ is the survival function of non-null $p$-value.  Assume,
\begin{itemize}
\item A1: For an appropriate choice of $\lambda$, $Q(\lambda)=P(p> \lambda|H=0)$ i.e the probability of non-null $p$-value being greater than $\lambda$ equals zero (Storey 2002).
\end{itemize} 
When parameter of interest under alternative hypothesis is substantially far from the same specified under null hypothesis or sample size is moderate to large, $p$-value tends to be smaller for consistent tests. Hence, even for moderate choice of $\lambda$, the probability of $p$-value under false null dominating $\lambda$, vanishes. This is a reasonable but crucial assumption in a sense that, violation of assumptions regarding the true value of the parameter of interest and sample size may not result in $Q(\lambda)=0$.
 Thus, applying A1 in (\ref{eq:Fbarl}) we get
 \begin{equation}\label{eq:pi0Fbarl}
 \pi_0=\frac{\bar{F}(\lambda)}{(1-\lambda)}.
 \end{equation} 
Let $p_1, p_2, ..., p_m$ be the $p$-values corresponding to the $m$ hypotheses tested or equivalently $m$ realizations on $p$. Denote $W(\lambda)=\sum_{i=1}^m I(p_i>\lambda)$ ($I$ denoting the indicator function) to be the number of $p$-values greater than $\lambda$. Putting the plug-in estimator of $\bar{F}(\lambda)$, i.e $W(\lambda)/m$ in (\ref{eq:pi0Fbarl}), an estimator for $\pi_0$ depending upon the choice of $\lambda$ may be suggested as
\begin{equation}\label{eq:storeyest}
  \hat{\pi}_0(\lambda)=\frac{W(\lambda)}{m\,(1-\lambda)}.  
\end{equation}
For a given dataset, two different choices of $\lambda$ would yield two different estimates and thus an optimum choice of $\lambda$ for a given dataset is necessary. For a subjectively chosen set with possible values of $\lambda\in \Lambda$, where $\Lambda=\{0, 0.05, 0.10, ..., 0.95\}$; a bootstrap routine is given in Storey (2002) and Storey et al. (2004) to approximate the best $\lambda$. Thus, Storey's bootstrap estimator is: $\hat{\pi}_0^B=\hat{\pi}_0(\lambda_{best})$. In Storey and Tibshirani (2003), natural cubic spline has been fitted to the $(\lambda, \hat{\pi}_0(\lambda))$ curve for smoothing and the evaluated value of the fit at $\lambda=1$ (as motivated in Corollary 1 of Storey (2002)) is taken as the final estimate which we denote by $\hat{\pi}_0^P$.\\ 

For a small choice of $\lambda$ in $\hat{\pi}_0(\lambda)$, the bias of the estimator is large while the variance is small. The situation is exactly opposite for large $\lambda$. It has been first noted by Jiang and Doerge (2008) and they  have suggested the use of multiple $\lambda$'s instead of a single best choice, in some sense. For the time being assume a fixed set $S_{\lambda}=\{(\lambda_1, \lambda_2, ..., \lambda_k): 0<\lambda_1< \lambda_2< ...< \lambda_k<1\}$ for a fixed $k$ and equal width given by $(\lambda_{i+1}-\lambda_i)$ for $i=1, 2, ..., k-1$ such that A1 holds. Then, the average estimate based approach suggests $\hat{\pi}_0^A=(1/k)\sum_{i=1}^k\,\hat{\pi}_0(\lambda_i)$ to be an appropriate estimator for $\pi_0$. The authors have also suggested a change-point based algorithm to select $S_\lambda$.
 
\subsection{Bias correction of Storey's estimator}    
Without the assumption A1, from (\ref{eq:Fbarl}) we get
\begin{equation}\label{pi0Fbarl2}
\pi_0=\frac{\bar{F}(\lambda)-Q(\lambda)}{(1-\lambda)-Q(\lambda)}
\end{equation} 
for fixed $\lambda$. Cheng et al. (2015) obtained (\ref{pi0Fbarl2}) from a somewhat different motivation. Substituting plug-in estimator of $\bar{F}(\lambda)$ has already been discussed in subsection 2.1. For estimating $Q(\lambda)$ following assumptions are necessary.
\begin{itemize}
\item A2: The availability of a common test for all the $m$ hypotheses. 
\item A3: The data-arrays used for each test are generated from a known parametric family.
\item A4: The closed form distribution of the test-statistics under the null are of a known family , enabling the calculation of the exact $p$-values.  
\item A5: The distribution of the non-null test-statistics and hence the non-null $p$-values are labelled by unknown effect sizes, which are different for each test. 
\end{itemize}
A2 is generally true for microarray experiments and is also appropriate for the present setup. Cheng et al. (2015) assumed normality for each expression level, such that A3 is valid. Time to events or its equivalent entities for each segment are assumed to be exponentially  distributed, thus satisfying A3. Under normality, the test-statistics for usual single-sample or two-sample tests for the mean are normal under null. In this work, the test-statistics for single-sample test related to the exponential rate parameter is a $\chi^2$ variate under null and for a two-sample problem the test-statistic is distributed as a $F$ variate. Thus, A4 also holds good. As mentioned earlier test for $H_i$ is performed by $T_i$ and we introduce the notation $\delta_i$ to denote the effect size of the corresponding test. Non-null distribution of $T_i$ and hence the non-null distribution of $p_i$ is to be labelled by $\delta_i$, $i=1, 2, ..., m$. Hung et al. (1997) have discussed properties of non-null $p$-values, where non-null distribution of the $p$-value for $Z$-test has been explored. For single sample and two sample $t$-test, similar discussion is available in section 3 of Cheng et al. (2015). We will discuss such properties of non-null $p$-value for single and two sample problems under exponential set-up in section 3. \\

Let $\mathcal{I}=\{1, 2, ..., m\}$. Also let $\mathcal{T}$ denote the set of indices corresponding to the originally true null hypotheses i.e, $\mathcal{T}=\{i\in \mathcal{I}: H_i=1\}$. Thus, the cardinality of $\mathcal{T}$, is $m_0$. Similarly denote the set of originally false null hypotheses by $\mathcal{F}$. Clearly, $\mathcal{F}=\mathcal{I}-\mathcal{T}$ with cardinality $m_1$. Each null $p$ has the  same distribution, uniform over $(0,1)$; while the distribution of non-null $p$-values are different but they belong to the same family. Let $f_1^\delta (p)$ denote the distribution of $p$ with effect size $\delta$. Then for all $i\in \mathcal{F}$, $Q_{\delta_i}(\lambda)=\int_\lambda^1 f_1^{\delta_i} (p)\, dp$, probability of $i$-th non-null $p$-value being greater than $\lambda$. Define, $Q^*(\lambda)=(1/m_1)\sum_{i\in \mathcal{F}} Q_{\delta_i}(\lambda)$, the average of non-null $p$-values greater than $\lambda$. By strong law of large number, $Q^*(\lambda)\to Q(\lambda)$, almost surely as $m_1\to\infty$. To estimate $Q(\lambda)$, individual $\delta_i$'s are estimated by $\hat{\delta}_i$, $i\in \mathcal{F}$. In fact, $\hat{\delta}_i$'s are strongly consistent for $\delta_i$ for each $i\in \mathcal{F}$. The estimation of $\delta$ under different testing problem is discussed in section 3. Each $Q_{\delta_i}(\lambda)$ is continuous in $\delta_i$ and thus, $Q_{\hat{\delta}_i }(\lambda)$ is strongly consistent for $Q_{\delta_i}(\lambda)$. Thus, a strongly consistent estimator for $Q(\lambda)$ is $\tilde{Q}(\lambda)=(1/m_1)\sum_{i\in \mathcal{F}} Q_{\hat{\delta}_i}(\lambda)$. In practice, $\mathcal{F}$ is unknown and hence $\tilde{Q}(\lambda)$ is unavailable.        Assume $\hat{Q}(\lambda)$ to be a dummy for $\tilde{Q}(\lambda)$ such that $\tilde{Q}(\lambda)\geq\hat{Q}(\lambda)$ with probability $ 1$. The computation of $\hat{Q}(\lambda)$ is discussed in detail in section 4. Substituting the plug-in estimators for $\bar{F}(\lambda)$ and $Q(\lambda)$ in (\ref{pi0Fbarl2}), we get $\tilde{\pi}_0^U(\lambda)$ (or $\hat{\pi}_0^U(\lambda)$), bias corrected estimator for $\pi_0$ with fixed choice of $\lambda$. We now address the issues related to reduction in bias and over-correction in the following result.\\

\noindent \textbf{Result 1}: With the set-up and notations introduced in subsection 2.2, for all $\lambda\in (0,1)$
\begin{eqnarray*}
&(a)\quad & \textrm{For}\,\, W(\lambda)/m\leq (1-\lambda),\,\, \tilde{\pi}^U_0(\lambda)\leq\hat{\pi}^U_0(\lambda)\leq\hat{\pi}_0(\lambda) .\\
&(b)\quad & \tilde{\pi}_0^U(\lambda)\,\to\,\pi_0,\,\textrm{almost surely} .
\end{eqnarray*}
\textbf{Proof}: Consider,
\begin{equation*}
    g(x)=\frac{a-x}{b-x}.
\end{equation*}
Note that, $g$ is non-increasing in $x$
 for $a\leq b$. Let $a=W(\lambda)/m$ and $b=(1-\lambda)$. Since, $0\leq\hat{Q}(\lambda)\leq\tilde{Q}(\lambda)$, $g(0)\geq g(\hat{Q}(\lambda))\geq g(\tilde{Q}(\lambda))$, which proves $(a)$. Now,
\begin{equation*}
    \tilde{\pi}^U_0(\lambda)=\frac{\frac{W(\lambda)}{m}-\tilde{Q}(\lambda)}{(1-\lambda)-\tilde{Q}(\lambda)}.
\end{equation*}
As $\hat{\delta}_i\xrightarrow{a.s} \delta_i\,\,\forall i=1, 2, ..., m$, $\tilde{Q}(\lambda)\xrightarrow{a.s}Q^*(\lambda)$. Thus, 
\begin{equation*}
    \tilde{\pi}_0^U(\lambda)\xrightarrow{a.s} \frac{\frac{W(\lambda)}{m}-Q^*(\lambda)}{(1-\lambda)-Q^*(\lambda)}=\tilde{\pi}_0,\,\,\, \textrm{say}.
\end{equation*}
Note that, $W(\lambda)/m\xrightarrow{a.s}\bar{F}(\lambda)$ as $m\to \infty$. As $m\to\infty$, $m_1\to\infty$ for $\pi_0\in (0, 1)$ and thus, $Q^*(\lambda)\xrightarrow{a.s} Q(\lambda)$ as $m\to\infty$. Hence $\tilde{\pi}_0\xrightarrow{a.s} (\bar{F}(\lambda)-Q(\lambda))/((1-\lambda)-Q(\lambda))=\pi_0$ as $m\to\infty$. $\blacksquare$\\

The situations $\hat{\pi}_0(\lambda)\leq 1$ and $\tilde{Q}(\lambda)\leq (1-\lambda)$ are quite usual in multi[ple testing set-up as the first one is a reasonable estimate of $\pi_0$ and $\tilde{Q}(\lambda)$ is a consistent estimate of $Q(\lambda)$, which is obviously less than $(1-\lambda)$. If these do not hold good, $\hat{\pi}^U_0(\lambda)$ lies outside the parameter space and then we take the estimate to be the nearest boundary point.\\

Result 1 combines claims written in section 2 and subsection 4.2 of Cheng et al. (2015). We have been able to prove Result 1 in a more direct way. Thus, the approach reduces conservative bias of Storey's primary estimator while refraining from over-correction.\\

$\Lambda=\{0.20, 0.25,..., 0.5\}$ is taken as in Jiang and Doerge (2008) for similar purpose (see subsection 2.2 in Cheng et al. (2015)) and we identify the following estimator as the bias and variance reduced estimator for $\pi_0$:
\begin{equation*}
\hat{\pi}_0^U=\frac{1}{\# \Lambda}\sum_{\lambda_j\in \Lambda} min\{1,max\{0,\hat{\pi}_0^U(\lambda_j)\}\} 
\end{equation*}
, where $\#\Lambda$ denotes cardinality of $\Lambda$.
\subsection{Estimator based on sum of all $p$-values}
Instead of taking $V$ degenerated at some fixed $\lambda$, assume $V\sim Uniform(0,1)$. Putting $G(v)=v$ for $v\in (0,1)$ in (\ref{eq:pgreaterV}) we get
\begin{equation}\label{eq:Ep}
P(p\geq V)=\int_0^1 p\,f(p)\,dp=E(p)=\frac{\pi_0}{2}+(1-\pi_0)\,e
\end{equation} 
since $p|H=1\sim Uniform(0,1)$. In (\ref{eq:Ep}) we use $e$ to denote expectation of non-null $p$-value: $e=E(p|H=0)$. From (\ref{eq:Ep}) we get
\begin{equation}\label{pi0Ep}
\pi_0=\frac{E(p)-e}{0.5-e}.
\end{equation}
To estimate $\pi_0$, both $E(p)$ and $e$ are to be estimated. $E(p)$ can be estimated by the mean of observed $p$-values: $\bar{p}=(1/m)\sum_{i=1}^m p_i$. Define $e^*=(1/m_1)\sum_{i\in \mathcal{F}}e_{\delta_i}$, which converges almost surely to $e$ as $m_1\to \infty$. $e$ can be estimated imitating the approach of estimating $Q(\lambda)$ with assumptions A2-A5. The corresponding estimator for $\pi_0$ has recently been introduced in Biswas (2019) and computation of $e_{\delta_i}=E(p_i|i\in \mathcal{F})$ has been demonstrated for single and two sample $t$-tests therein. Since each $e_{\delta_i}$ is bounded and continuous in $\delta_i$, following the discussion in subsection 2.2, a strongly consistent estimator for $e$ is $\tilde{e}=(1/m_1)\sum_{i\in \mathcal{F}}e_{\hat{\delta}_i}$, which cannot be realized in practice for obvious reason mentioned earlier and hence $\tilde{\pi}_0^E=(\bar{p}-\tilde{e})/(0.5-\tilde{e})$ cannot be implemented. For $\hat{e}$ being a dummy of $\tilde{e}$ with $\hat{e}\leq \tilde{e}$ almost surely, an workable estimator for $\pi_0$ is $\hat{\pi}_0^E=(\bar{p}-\hat{e})/(0.5-\hat{e})$.\\

\noindent \textbf{Result 2}: With the set-up and notations introduced in subsection 2.3
\begin{eqnarray*}
&(a)\quad & \textrm{For}\,\, \bar{p}\leq 0.5,\,\, \tilde{\pi}_0^E\leq \hat{\pi}_0^E.\\
&(b)\quad & \tilde{\pi}_0^E\,\to\,\pi_0,\,\textrm{almost surely}.
\end{eqnarray*}  
\textbf{Proof}: We consider $g$ as in Result 1 and assume $a=\bar{p}$, $b=0.5$. Since, $\hat{e}\leq\tilde{e}$, $g(\hat{e})\geq g(\tilde{e})$, which proves $(a)$. Here, 
\begin{equation*}
    \tilde{\pi}^E_0=\frac{\bar{p}-\tilde{e}}{0.5-\tilde{e}}.
\end{equation*}
As $\hat{\delta}_i\xrightarrow{a.s}\delta_i\,\,\, \forall i=1, 2, ..., m$, $\tilde{e}\xrightarrow{a.s} e^*$. Thus,
\begin{equation*}
    \tilde{\pi}_0^E\xrightarrow{a.s}\frac{\bar{p}-e^*}{0.5-e^*}=\tilde{\pi}_0,\,\,\textrm{say}.
\end{equation*}
Note that, $\bar{p}\xrightarrow{a.s} E(p)$ as $m\to\infty$. As $m\to \infty$, $m_1\to \infty$ and thus, $e^*\xrightarrow{a.s} e$ as $m\to\infty$. Hence, $\tilde{\pi}_0\xrightarrow{a.s}(E(p)-e)/(0.5-e)=\pi_0$ as $m\to \infty$.$\blacksquare$\\

The situations $\bar{p}\leq 0.5$ and $\tilde{e}\leq 0.5$ are very natural in multiple testing set-up as $\bar{p}$ is consistent for $E(p)$, which is less than $0.5$ and similarly $\tilde{e}$ is consistent for $e$ which is also less than $0.5$. If theses do not hold good, $\hat{\pi}^E_0$ lies outside the parameter space and then we take the estimator for $\pi_0$ as
\begin{equation*}
    \hat{\pi}^E_0=min\left\{1, max \left\{0, \frac{\bar{p}-\hat{e}}{0.5-\hat{e}}\right\}\right\}.
\end{equation*}

Both the model based bias corrected estimators are shown to have conservative bias for estimating $\pi_0$. In Cheng et al. (2015), $\hat{\pi}_0^U$ has been shown to outperform the robust estimators under reasonable model assumption whereas, under similar situation $\hat{\pi}_0^E$ outperforms it in terms of mean square error, as empirically studied in Biswas (2019) through extensive simulation study. Note that, both the estimators use an initial estimator for $\pi_0$ but the computation of $\hat{\pi}_0^E$ does not require flexible threshold tuning parameters owing to the fact that it uses all the $p$-values. To rule out the possibility of estimates taking value outside the parameter space under very unusual situation, $\hat{\pi}_0^E$ is taken to be equal to the nearest boundary point when it lies outside the parameter space.

\section{Properties of non-null $p$-values}
To implement the bias corrected estimators $\hat{\pi}_0^U$ and $\hat{\pi}_0^E$, appropriate estimates of the unknown quantities $Q(\lambda)$ and $e$ are needed. To get explicit expressions for these quantities, we need to have the probability density functions $f_1^{\delta_i}(p)$ (for notational convenience we write this to be $f_{\delta_i}(p)$ henceforth) for each non-null $p$-value with effect size $\delta_i$, $i\in \mathcal{F}$. The subscript $i$ in effect sizes are not specified in this section for ease of notation. Thus, for different testing scenarios we determine the probability density function $f_\delta(p)$, then $Q_\delta(\lambda)$ by integrating $f_\delta(p)$ from $\lambda$ to $1$ and finally obtain $e_\delta$ through the following results. As discussed in subsection 2.2, $Q_{\delta}(\lambda)$ for fixed $\lambda$ and $e_\delta$ are continuous in $\delta$ under each of the testing problems considered here.\\

\noindent\textbf{Result 3}: Assume $X_1, X_2, ..., X_n$ be a random sample of size $n$ from an exponential distribution with mean $\theta$. Consider the following testing problem:
\begin{equation}\label{eq:hyp1}
H_0: \theta=1\quad \quad\textrm{versus}\quad\quad H_0:\theta>1.
\end{equation}
For the corresponding likelihood ratio test
\begin{eqnarray*}
&(a)\quad &\delta=\theta\quad \textrm{and thus}\quad \hat{\delta}=min\{1,\bar{X}\}\\
&(b)\quad & f_\delta(p)=\frac{\frac{1}{\delta}\,f_{\chi_{2n}^2}(\frac{1}{\delta}\,\chi_{p,2n}^2)}{f_{\chi_{2n}^2}(\chi_{p,2n}^2)}\quad \textrm{for}\quad 0<p<1\\
&(c)\quad & Q_\delta(\lambda)=F_{\chi_{2n}^2}\left(\frac{1}{\delta}\, \chi_{\lambda,2n}^2\right)\quad \textrm{for}\quad 0<\lambda <1\\
&(d)\quad & e_\delta=E_{X\sim\chi_{2n}^2}\left[1-F_{\chi_{2n}^2}\left(\frac{X}{\delta}\right)\right].
\end{eqnarray*}
Here $f_{\chi_{\nu}^2}$, $F_{\chi_{\nu}^2}$ and $\chi_{p,\nu}^2$ denote the probability density function, the distribution function and the upper-$p$ point of chi-square distribution with $\nu$ degrees of freedom, respectively.\\

\noindent\textbf{Proof}: The likelihood ratio test corresponding to the hypothesis in (\ref{eq:hyp1}) uses the test-statistic $T=2\sum_{i=1}^n X_i \sim \theta\chi_{2n}^2$. Effect size of the test $\delta=\theta$. As $E(T)=2n\delta$, an unbiased estimator of $\delta$ is $\hat{\delta}=T/2n=\bar{X}$, the sample mean. \\

As we reject $H_0$ for larger observed value of $T$, the corresponding $p$-value is defined as $p=P_{H_0}(\chi_{2n}^2>T)=1-F_{\chi_{2n}^2}(T)$ since under $H_0$, $T\sim \chi_{2n}^2$. Therefore, $p\sim Uniform(0,1)$, under $H_0$. Under $H_1$, $T\sim \delta\chi_{2n}^2$ and therefore the density function of $T$, labelled by $\delta$ is 
\begin{equation}\label{eq:tFinv1}
f_{\delta}(t)=\frac{1}{\delta}f_{\chi_{2n}^2}\left(\frac{t}{\delta}\right)\quad \textrm{for}\quad t>0.
\end{equation}    
From the relation between $T$ and $p$, $t=F_{\chi_{2n}^2}^{-1}(1-p)=\chi_{p,2n}^2$. The corresponding absolute Jacobian of transformation is $f_{\chi_{2n}^2}(\chi_{p,2n}^2)$. Thus from (\ref{eq:tFinv1}), the density function of $p$ labelled by $\delta$ is
\begin{equation}\label{eq:fp1}
f_\delta(p)=\frac{\frac{1}{\delta}f_{\chi_{2n}^2}(\chi_{p,2n}^2)}{f_{\chi_{2n}^2}(\chi_{p,2n}^2)}\quad\textrm{for}\quad 0<p<1.
\end{equation} 

For $\lambda\in (0,1)$ upper tail probability labelled by $\delta$, using (\ref{eq:fp1}) in expression of $Q_\delta(\lambda)$ we get
\begin{eqnarray*}
Q_\delta(\lambda)&=&\int_{\lambda}^1 \frac{\frac{1}{\delta}f_{\chi_{2n}^2}(\frac{1}{\delta}\chi_{p,2n}^2)}{f_{\chi_{2n}^2}(\chi_{p,2n}^2)}\,dp\,\, =\,\,I,\,\, \textrm{say}.\\
\end{eqnarray*}
By change of variable from $p$ to $v$ such that $v=(1/\delta)\chi_{p,2n}^2$ we get
\begin{equation*}
    I=F_{\chi_{2n}^2}(\frac{1}{\delta}\chi_{\lambda,2n}^2).
\end{equation*}
which proves the result in $(c)$. For an explicit expression for expected $p$-value under the false null, we apply \\
\begin{equation*}
e_\delta=\int_0^1 F_{\chi_{2n}^2}\left(\frac{1}{\delta}\chi_{p,2n}^2\right)\,dp
\end{equation*}
By change of variable from $p$ to $v$ such that $v=\chi_{p,2n}^2$ we get
\begin{eqnarray*}
e_\delta &=& \int_0^\infty F_{\chi_{2n}^2}\left(\frac{v}{\delta}\right)\, f_{\chi_{2n}^2}(v)\,dv\\
&=& E_{X\sim \chi_{2n}^2} \left[F_{\chi_{2n}^2}\left(\frac{X}{\delta}\right)\right].\qquad\blacksquare
\end{eqnarray*}

\noindent \textbf{Result 4}: For $X_1, X_2, ..., X_n$ to be a random sample from exponential distribution with mean $\theta$, consider the testing problem:
\begin{equation}\label{eq:hyp2}
H_0:\theta=1\quad\textrm{versus}\quad H_1:\theta\neq 1.
\end{equation} 
For the corresponding likelihood ratio test
\begin{eqnarray*}
&(a)\quad & \delta=\theta\quad\textrm{and thus}\quad \hat{\delta}=\bar{X}\\
&(b)\quad & Q_\delta(\lambda)=F_{\chi_{2n}^2}\left[\frac{1}{\delta}\chi_{\frac{\lambda}{2},2n}^2\right]-F_{\chi_{2n}^2}\left[\frac{1}{\delta}\chi_{1-\frac{\lambda}{2},2n}^2\right]\quad\textrm{for}\quad 0<\lambda<1\\
&(c)\quad & e_\delta=E_{X\sim \chi_{2n}^2(0,\mu)}\left[F_{\chi_{2n}^2}\left(\frac{X}{\delta}\right)\right]-E_{X\sim \chi_{2n}^2(\mu,\infty)}\left[F_{\chi_{2n}^2}\left(\frac{X}{\delta}\right)\right].
\end{eqnarray*}
The notations used for stating Result 3 also remain relevant here. In addition to that, $\chi_{\nu}^2(a,b)$ denotes the truncated chi-squared distribution with degree of freedom $\nu$ and region of truncation being $(a,b)$. Here $\mu$ denotes the median of $\chi_{2n}^2$ distribution.  \hspace{5.8cm}$\blacksquare$\\

\noindent \textbf{Proof}: The corresponding likelihood ratio test uses the same test-statistic as in Result 3 and thus part $(a)$ of Result 4 follows directly from part $(a)$ of Result 3. For the next part, it should be noted that due to two-sided alternative hypotheses, the corresponding $p$-value is defined through $T$, where
\begin{eqnarray*}
p&=&2\,min\left\{P(\chi_{2n}^2>T),P(\chi_{2n}^2<T)\right\}\\
&=&2\,min(p^*,1-p^*),\, \textrm{say}.
\end{eqnarray*}
Here $p^*=P(\chi_{2n}^2>T)$ is the $p$-value defined for the testing problem in (\ref{eq:hyp1}). Thus from part $(b)$ of Result 3 we have
\begin{equation}\label{eq:fpstar}
f_\delta(p^*)=\frac{\frac{1}{\delta}f_{\chi_{2n}^2}(\chi_{p^*,2n}^2)}{f_{\chi_{2n}^2}(\chi_{p^*,2n}^2)}.
\end{equation}
Now for any $\lambda\in (0,1)$,
\begin{eqnarray*}
Q_\delta(\lambda)&=&P(p>\lambda)\\
&=&P\left(\frac{\lambda}{2}<p^*<1-\frac{\lambda}{2}\right)\\
&=& \int_{\frac{\lambda}{2}}^{1-\frac{\lambda}{2}}\frac{\frac{1}{\delta}f_{\chi_{2n}^2}(\chi_{p^*,2n}^2)}{f_{\chi_{2n}^2}(\chi_{p^*,2n}^2)}\,dp^*\quad\textrm{[using (\ref{eq:fpstar})]}\\
&=&\int_{\frac{1}{\delta}\chi_{1-\frac{\lambda}{2},2n}^2}^{\frac{1}{\delta}\chi_{\frac{\lambda}{2},2n}^2} f_{\chi_{2n}^2}(v)\,dv \quad\textrm{by taking}\quad v=\frac{1}{\delta}\chi_{p^*,2n}^2,
\end{eqnarray*}
which proves the result in $(b)$. 
\begin{eqnarray*}
e_\delta &=& \int_0^1 Q_\delta(\lambda)\, d\lambda\\
&=& \int_0^1 F_{\chi_{2n}^2}\left[\frac{1}{\delta}\chi_{\frac{p}{2},2n}^2\right]\,dp-\int_0^1 F_{\chi_{2n}^2}\left[\frac{1}{\delta}\chi_{1-\frac{p}{2},2n}^2\right]\,dp\\
&=& I_1-I_2,\,\textrm{say}.
\end{eqnarray*} 
Now, we consider the problem of evaluating the integral $I_1$. By change of variable from $p$ to $v$ such that $v=\chi_{\frac{p}{2},2n}^2$, we get
\begin{eqnarray*}
I_1&=& 2\int_{\chi_{\frac{1}{2},2n}^2}^\infty F_{\chi_{2n}^2}\left(\frac{v}{\delta}\right)\,f_{\chi_{2n}^2}(v)\,dv\\
&=& E_{X\sim \chi_{2n}^2(0,\mu)}\left[F_{\chi_{2n}^2}\left(\frac{X}{\delta}\right)\right]
\end{eqnarray*}
Following the same steps for evaluating $I_1$, $I_2$ can also be evaluated and thus the result in $(c)$. $\blacksquare$\\

\noindent \textbf{Result 5}: Let $X_1, X_2, ..., X_{n_1}$ and $Y_1, Y_2, ..., Y_{n_2}$ be two random samples of size $n_1$ and $n_2$ respectively from exponential distribution with mean $\theta_1$ and $\theta_2$. Consider the testing problem
\begin{equation}\label{eq:hyp3}
H_0:\theta_2=\theta_1\quad\textrm{versus}\quad H_1:\theta_2>\theta_1.
\end{equation}
For the corresponding likelihood ratio test, we have the following.
\begin{eqnarray*}
&(a)\quad & \delta=\frac{\theta_2}{\theta_1}\quad \textrm{and}\quad \hat{\delta}= min\left\{1,\frac{n_1-1}{n_1}\,\frac{\sum_{i=1}^{n_2}Y_i}{\sum_{i=1}^{n_1}X_i}\right\}\\
&(b)\quad & f_\delta(p)=\frac{\frac{1}{\delta}f_{F_{2n_2,2n_1}\left(\frac{1}{\delta}F_{p,2n_2,2n_1}\right)}}{f_{F_{2n_2,2n_1}\left(F_{p,2n_2,2n_1}\right)}}\quad\textrm{for}\quad 0<p<1\\
&(c)\quad & Q_\delta(\lambda)=F_{F_{2n_2,2n_1}}\left(\frac{1}{\delta}F_{\lambda,2n_2,2n_1}\right)\quad\textrm{for}\quad 0<\lambda<1\\
&(d)\quad & e_\delta =E_{X\sim F_{2n_2,2n_1}}\left[F_{F_{2n_2,2n_1}}\left(\frac{X}{\delta}\right)\right].
\end{eqnarray*}
Here $f_{F_{\nu_1,\nu_2}}$, $F_{F_{\nu_1,\nu_2}}$ and $F_{p,\nu_1,\nu_2}$ respectively denote the probability density function, the distribution function and the upper-$p$ point of $F$ distribution with $\nu_1$ and $\nu_2$ degrees of freedom.\\

\noindent \textbf{Proof}: The likelihood test corresponding to the hypothesis in (\ref{eq:hyp3}) uses the test-statistic $T=\sum_{i=1}^{n_2} Y_i/\sum_{i=1}^{n_1} X_i\sim (\theta_2/\theta_1) F_{2n_2,2n_1}$. Effect size of the test $\delta=\theta_2/\theta_1$. As $E(T)=\delta[n_1/(n_1-1)]$, an unbiased estimator of $\delta$ is $\hat{\delta}=[(n_1-1)/n_1]T$ and thus the result in $(a)$ follows. Rest of the proof follows from the proof of Result 3 with obvious changes.$\blacksquare$\\

\noindent \textbf{Result 6}: Consider $X_1, X_2, ..., X_{n_1}$ and $Y_1, Y_2, ..., Y_{n_2}$ to be independent random samples of size $n_1$ and $n_2$ from exponential distributions with mean $\theta_1$ and $\theta_2$, respectively. Consider the testing problem
\begin{equation}\label{eq:hyp4}
H_0:\theta_2=\theta_1\quad\textrm{versus}\quad H_1:\theta_2\neq\theta_1.
\end{equation}
For the corresponding likelihood ratio test, we have 
\begin{eqnarray*}
&(a)\quad & \delta=\frac{\theta_2}{\theta_1}\quad \textrm{and} \quad \hat{\delta}=\frac{n_1-1}{n_1}\,\frac{\sum_{i=1}^{n_2}Y_i}{\sum_{i=1}^{n_1}X_i}\\
&(b)\quad & Q_\delta(\lambda)=F_{F_{2n_2,2n_1}}\left[\frac{1}{\delta}F_{\frac{\lambda}{2},2n_2,2n_1}\right]-F_{F_{2n_2,2n_1}}\left[\frac{1}{\delta}F_{1-\frac{\lambda}{2},2n_2,2n_1}\right]\\
&(c)\quad & e_\delta=E_{X\sim F_{2n_2,2n_1}(0,\mu)}\left[F_{F_{2n_2,2n_1}}\left(\frac{X}{\delta}\right)\right]-E_{X\sim F_{2n_2,2n_1}(\mu,\infty)}\left[F_{F_{2n_2,2n_1}}\left(\frac{X}{\delta}\right)\right].
\end{eqnarray*} 
The notations used for Result 5 also remain relevant here. In addition to that, $F_{\nu_1,\nu_2}(a,b)$ denotes the truncated $F$-distribution with degrees of freedom $\nu_1$, $\nu_2$ and region of truncation $(a,b)$. Here $\mu$ denotes the median of $F_{2n_2,2n_2}$ distribution.\\

\noindent \textbf{Proof}: For the testing problem in (\ref{eq:hyp4}), the likelihood ratio test uses the same test-statistic as in Result 5. Since the critical region is two-sided, the corresponding $p$-value is similarly defined as in Result 4. One can follow the steps elaborated through the proof of Result 4 and use Result 5 to easily prove Result 6. $\blacksquare$

\section{Algorithms}
Algorithm for computing $\hat{\pi}_0^U$ under normal model assumption is given in in Cheng et al. (2015) and for $\hat{\pi}_0^E$ under same set-up, see  Biswas (2019). First, we reframe the algorithms under current set-up to maintain readability and for making the proposed estimation methods readily available to the practitioners. For the sake of brevity we only consider the testing problem in Result 4 and use the corresponding non-null $p$-value properties here. For all the four situations discussed here, the following algorithms can be implemented with obvious modifications.
\begin{center}
\textbf{Algorithm 1}\\
(For computing $\hat{\pi}_0^U$)
\end{center}
\begin{itemize}
\item For all $i=1, 2, ..., m$, estimate $\delta_i$ by $\hat{\delta}_i=\bar{X}_i$.
\item For all $i=1, 2, ..., m$ and for each $\lambda_j\in \Lambda$;  estimate the upper tail probability $Q_{\delta_i}(\lambda_j)$ by $Q_{\hat{\delta}_i}(\lambda_j)$ given by
\begin{equation*}
Q_{\hat{\delta}_i}(\lambda_j)=F_{\chi_{2n_i}^2}\left[\frac{1}{\hat{\delta}_i}\chi_{\frac{\lambda_j}{2},2n_i}^2\right]-F_{\chi_{2n_i}^2}\left[\frac{1}{\hat{\delta}_i}\chi_{1-\frac{\lambda_j}{2},2n_i}^2\right].
\end{equation*} 
where, $n_i$ denotes the available sample size for testing $i$-th hypothesis.
\item Using an available estimator of $\pi_0$ as initial estimator $\hat{\pi}_0^I$, calculate $d=[m\times (1-\hat{\pi}_0^I)]$, where $[\,\,]$ denotes the usual box function. Arrange $Q_{\hat{\delta}_i}(\lambda_j)$'s in increasing order and denote the $i$-th quantity in the list as $\hat{Q}_{(i)}(\lambda_j)$. Thus a conservative estimator for $Q(\lambda_j)$ is
\begin{equation*}
\hat{Q}(\lambda_j)=\frac{1}{d}\sum_{i=1}^d \hat{Q}_{(i)}(\lambda_j).
\end{equation*}
\item Given $\hat{Q}(\lambda_j)\,\,\,\forall \,\,\,\lambda_j\in \Lambda$, calculate 
\begin{equation*}
\hat{\pi}_0^U=\frac{1}{\#\Lambda}\sum_{\lambda_j\in\Lambda} min\left\{1,max\left\{0,\frac{\frac{W(\lambda_j)}{m}-\hat{Q}(\lambda_j)}{(1-\lambda_j)-\hat{Q}(\lambda_j)}\right\}\right\}.
\end{equation*}
\end{itemize}

\begin{center}
\textbf{Algorithm 2}\\
(For computing $\hat{\pi}_0^E$)
\end{center}
\begin{itemize}
\item For all $i=1, 2, ..., m$, estimate $\delta_i$ by $\hat{\delta}_i=\bar{X}_i$.
\item For all $i=1, 2, ..., m$, estimate the mean of non-null $p$-value $e_{\delta_i}$ by
\begin{equation*}
\hat{e}_{\delta_i}= E_{X\sim \chi_{2n_i}^2(0,\mu_i)}\left[F_{\chi_{2n_i}^2}\left(\frac{X}{\hat{\delta}_i}\right)\right]-E_{X\sim \chi_{2n_i}^2(\mu_i,\infty)}\left[F_{\chi_{2n_i}^2}\left(\frac{X}{\hat{\delta}_i}\right)\right].
\end{equation*}
where, $n_i$ denotes the available sample size for testing $i$-th hypothesis and $\mu_i$ denotes the median of $\chi_{2n_i}^2$ distribution.
\item Using an initial estimator of $\pi_0$ as initial estimator $\hat{\pi}_0^I$, calculate $d=[m\times(1-\hat{\pi}_0^I)]$, as before. Arrange $\hat{e}_{\delta_i}$'s in increasing order and denote the $i$-th quantity in the list as $\hat{e}_{(i)}$. Thus a conservative estimator for $e$ is 
\begin{equation*}
\hat{e}=\frac{1}{d}\sum_{i=1}^d  \hat{e}_{(i)}.
\end{equation*} 
\item Given $\hat{e}$, calculate
\begin{equation*}
\hat{\pi}_0^E=min\left\{1,max\left\{0,\frac{\bar{p}-\hat{e}}{0.5-\hat{e}}\right\}\right\}.
\end{equation*}
\end{itemize}

\noindent \textbf{Note 1}: The role of $\hat{\pi}_0^I$ is important in obtaining $\hat{Q}(\lambda)$ and $\hat{e}$. For $\hat{\pi}_0^I\geq \pi_0$, observe that $m_1\geq d$. Clearly, $\hat{Q}(\lambda)\leq \tilde{Q}(\lambda)$ and $\hat{e}\leq \tilde{e}$. Thus, $\hat{\pi}_0^U\geq \tilde{\pi}_0^U$ and $\hat{\pi}_0^E\geq \tilde{\pi}_0^E$.\\

\noindent \textbf{Note 2}: For implementation of both the algorithms, we choose Storey's bootstrap estimator, given by $\hat{\pi}_0^B$ as the initial estimator. This choice seems reasonable albeit being non-universal and further research on this is warranted. The algorithms could also be implemented with other choices of $\hat{\pi}_0^I$. The performance analysis of the bias corrected estimators under the current set-up requires extensive simulation study, starting with different choices of the initial estimator. In fact, the algorithms could in principle be done several times, each time with the estimate of $\pi_0$ from the previous iteration. Obviously, this technique will become computation intensive for all practical purposes. We refrain from addressing these issues, as they are beyond the scope of the current study.\\

It has already been mentioned in section 1 that Benjamini-Hochberg procedure for controlling the false discovery rate is conservative. To understand this, we briefly discuss FDR and the algorithm for controlling it at a prefixed level $q\in (0,1)$. While testing $m$ hypotheses simultaneously, let $R$ be the total number of rejected hypotheses by application of certain multiple testing algorithm. From the entire set of rejected hypotheses, some hypotheses may be originally true. These are categorized as false discovery and let $V$ denote the total number of such false discoveries. Then the false discovery proportion (FDP) is defined as 
\begin{equation*}
FDP=\begin{cases}
\frac{V}{R}\quad \textrm{if}\quad R>0\\
0\quad \textrm{if}\quad R=0.
\end{cases}
\end{equation*}    
Note that, prior to the application of any algorithm both $V$ and $R$ are random variables and the expected value of FDP is termed as $ \,\,False \,\,Discovery\,\, Rate$ (FDR). Let $p_{(1)}\leq p_{(2)}\leq...\leq p_{(m)}$ be the ordered sequence of the available $p$-values. Benjamini-Hochberg procedure identifies the largest $k$ such that $p_{(k)}\leq \left(k/m\right)q$ and rejects all hypotheses with corresponding $p$-value less than $p_{(k)}$ along with the hypothesis with $p$-value $p_{(k)}$. This procedure is conservative, as the  implementation of the same ensures $FDR=\pi_0\,q$ where, $\pi_0=m_0/m$. To overcome this shortcoming, Craiu and Sun (2008) worked with the following adaptive Benjamini-Hochberg procedure which uses an approximation of $\pi_0$. 
\begin{center}
\textbf{Algorithm 3}\\
(Implementing adaptive BH procedure to control FDR at level $q$)
\end{center}    
\begin{itemize}
\item Let the $p$-value corresponding to the problem of testing $H_i$ be $p_i$ for $i=1, 2, ..., m$. Arrange the available $p$-values in increasing order: $p_{(1)}, p_{(2)}, ..., p_{(m)}$. Denote the corresponding hypotheses by $H_{(i)}:\quad i=1, 2, ..., m$.
\item Given the data-set, estimate $\pi_0$. Let it be $\hat{\pi}_0$.
\item Compute the adjusted $p$-value corresponding to $p_{(i)}$:
\begin{equation*}
adj.p_{(i)}=min\left\{\hat{\pi}_0\frac{m\,p_{(j)}}{j}:j\geq i\right\}\quad \forall\quad i=1, 2, ..., m.
\end{equation*}
\item For all $i=1, 2, ..., m$, reject $H_{(i)}$if $adj.p_{(i)}\leq q$.
\end{itemize}
Both adaptive BH procedure and Storey's q-value approach are justified to be equivalent in Craiu and Sun (2008). They have also emphasized that both the approaches require a good approximation of $\pi_0$. Less conservative estimators for $\pi_0$ are in demand since closer approximation of $\pi_0$ will bring superiority in the adaptive procedure by increasing the number of rejections while controlling FDR at level $q$, as evident from Algorithm 3. In numerical study, we use \begin{slshape} adjust.p( )\end{slshape} function (available in \begin{slshape}cp4p\end{slshape} library from Bioconductor) by Gianetto et al. (2019) for obtaining the adjusted $p$-values.

\section{Simulation study}
We have conducted an extensive simulation study to investigate the performance of the bias corrected estimators under different settings. The well-known and established estimators apart from the proposed $\hat{\pi}_0^U$ and $\hat{\pi}_0^E$, considered for performance comparison are listed below. 
\begin{eqnarray*}
\hat{\pi}_0^B &:& \textrm{Storey's bootstrap estimator (discussed in subsection 2.1)}\\
\hat{\pi}_0^L &:& \textrm{Convest estimator (Langaas et al. 2005)}\\
\hat{\pi}_0^A &:& \textrm{Jiang and Doerge's average estimator (discussed in subsection 2.1)}\\
\hat{\pi}_0^P &:& \textrm{ Natural cubic spline smoothing based estimator (discussed in subsection 2.1)}\\
\hat{\pi}_0^H &:& \textrm{Histogram based estimator (Nettleton et al. 2006)}\\
\hat{\pi}_0^D &:& \textrm{A robust estimator of $\pi_0$ (Pounds and Cheng 2006)}\\
\hat{\pi}_0^S &:& \textrm{Sliding linear model based estimator (Wang et al. 2011)}.
\end{eqnarray*}  
\subsection{Simulation setting}
We imitate a segmented time to event dataset to generate artificial datasets. For this purpose we choose $m=100, 500, 1000$ segments and for each segment $n=30, 50$ available observations. For fixed $\pi_0=0.1, 0.2, ..., 0.9$, calculate $m_0=[m\,\pi_0]$ and take $m_1=m-m_0$. We set the mean failure time under null $\theta_0$ as `unity'. For $m_0$ randomly chosen numbers from $\mathcal{I}=\{1, 2, ..., m\}$ we fix $\theta=\theta_0=1$ and for the remaining cells in the array $\theta$, we generate values through some stochastic mechanism ensuring that they are not equal to $\theta_0$. We discuss this through the following two settings, although any other convenient settings might be used for the purpose.
\begin{itemize}
\item Uniform setting: For segments with better average lifetime $\theta\sim \textrm{Uniform}(1,1.5)$ and for segments with poor average lifetime $\theta\sim \textrm{Uniform}(0.5,1)$.

\item Exponential setting: For segments with better average lifetime $\theta$ is generated from truncated exponential distribution (mean=1) with support $(1,1.5)$ and for segments with poor average lifetime, the same is done with support $(0.5,1)$. 

\end{itemize}
Performance of an estimator for $\pi_0$ may vary with the proportion of better (or poor) non-null mean lifetimes. Thus we also replicate the experiment in each setting with different proportions of left \& right non-null $\theta$-values given respectively by the allocation proportions 75\% \& 25\%, 50\% \& 50\% and 25\% \& 75\%. \\

After generating the array of parameter $\theta$, we generate a sample of size $n$ from the exponential distribution with mean $\theta_i$ for all $i=1, 2, ..., m$. Thus the data matrix of order $m\times n$ is generated where each row correspond to observations from a particular segment and out of them $m_0$ (fixed by the choice of $\pi_0$) segments originally have mean lifetime equal to $\theta_0=1$. From each row of the data matrix we obtain $p$-value by applying appropriate test and construct a $p$-value array of length $m$ to compute the bias corrected estimators from Algorithm 1 and Algorithm 2. The other estimators are computed using \begin{slshape}estim.pi0\end{slshape} R-function (available in \begin{slshape}cp4p\end{slshape} library). Algorithm 3 also uses this array of $p$-values and an estimate of $\pi_0$ to identify the significantly different segments with control over FDR.    
\subsection{Simulation results}
Under different set-ups mentioned in subsection 5.1, each experiment is repeated $N=1000$ times and the estimators are compared through
\begin{eqnarray*}
MSE(\hat{\pi}_0)&=&\frac{1}{N}\sum_{i=1}^N (\hat{\pi}_{0i}-\pi_0)^2\quad \quad \textrm{and}\\
Bias(\hat{\pi}_0)&=&\frac{1}{N}\sum_{i=1}^N (\hat{\pi}_{0i}-\pi_0).
\end{eqnarray*}
MSE and bias of each estimators under uniform setting with $50\%-50\%$ allocation are reported in Table 1 and Table 2, respectively. From Table 1, it can be pointed out that the bias corrected estimators beat other estimators over a significant region of the parameter space (for $\pi_0\in (0,0.6)$) while $\hat{\pi}_0^U$ performs slightly better than $\hat{\pi}_0^E$. Thus, their performance may be considered to be approximately equivalent. Thus using the bias corrected estimators for small to moderate values of $\pi_0$ brings significant improvement while for larger values of $\pi_0$ it remains a viable alternative. From Table 2, similar comments may be made. Additionally, we point out that $\hat{\pi}_0^U$ really reduces the bias of $\hat{\pi}_0^B$ for a significant portion of the parameter space. As expected, the mean squared error for different estimators increase with increasing $m/n$ ratio, while relative performance of the proposed bias corrected estimators gets better when the same ratio increases. Observations under other allocation proportions of non-null $\theta$-values are almost similar and we refrain from reporting the numerical results to keep it concise. However, the gain from improved estimation of $\pi_0$ needs to be elaborated. Precise estimation of $\pi_0$ is used to apply adaptive algorithm for identifying significant segments as mentioned in section 4. We construct adaptive algorithms through the Algorithm 3 with different estimators for $\pi_0$ and compute the proportion of rejections for each which may be considered as the power of a particular algorithm in multiple testing set-up. In Figure 1 we compare power of non-adaptive Benjamini-Hochberg algorithm with its adaptive version using $\hat{\pi}_0$ as function of $\pi_0$ and observe that for lower to moderate values of $\pi_0$, the adaptive version results in substantial gain in power, while marginal gain is observed for larger values of $\pi_0$. Similar comparison for adaptive algorithm with different estimates of $\pi_0$ is reported in Table 3. From Table 3, it is evident that the bias corrected estimators outperform the others for lower to moderate values of $\pi_0$, where it really matters as pointed out from Figure 1.       

\begin{singlespacing}

\begin{small}
\begin{longtable}{cccccccccc}
\caption{MSE of estimators under different simulation settings}\\
\endfirsthead
\hline
\\
\multicolumn{10}{c}{m=100     n=50}\\
\\
\hline
$\pi_0$&$\hat{\pi}_0^B$&$\hat{\pi}_0^L$&$\hat{\pi}_0^A$&$\hat{\pi}_0^P$&$\hat{\pi}_0^H$&$\hat{\pi}_0^D$&$\hat{\pi}_0^S$&$\hat{\pi}_0^U$&$\hat{\pi}_0^E$\\
\hline
0.1&0.10995&0.11170&0.15401&0.21527&0.15302&0.17619&0.22819&\textbf{0.08046}&0.08484\\
0.2&0.08539&0.08641&0.12463&0.17941&0.12251&0.13958&0.17753&\textbf{0.06559}&0.06880\\
0.3&0.06782&0.06877&0.10510&0.16863&0.09821&0.11078&0.13427&\textbf{0.05506}&0.05762\\
0.4&0.04932&0.04902&0.08092&0.13429&0.07220&0.08275&0.09819&\textbf{0.04288}&0.04478\\
0.5&0.03465&0.03452&0.06114&0.10512&0.05136&0.05929&0.06630&\textbf{0.03188}&0.03322\\
0.6&0.02391&\textbf{0.02269}&0.04486&0.07593&0.03504&0.04003&0.04138&0.02345&0.02432\\
0.7&0.01600&\textbf{0.01330}&0.02905&0.04868&0.01999&0.02453&0.02311&0.01623&0.01640\\
0.8&0.01073&\textbf{0.00756}&0.01635&0.02949&0.00905&0.01265&0.00978&0.01015&0.01051\\
0.9&0.00994&0.00516&0.00642&0.01750&\textbf{0.00353}&0.00449&0.00947&0.00589&0.00617\\
\hline
\\
\multicolumn{10}{c}{m=100    n=30}\\
\\
\hline
$\pi_0$&$\hat{\pi}_0^B$&$\hat{\pi}_0^L$&$\hat{\pi}_0^A$&$\hat{\pi}_0^P$&$\hat{\pi}_0^H$&$\hat{\pi}_0^D$&$\hat{\pi}_0^S$&$\hat{\pi}_0^U$&$\hat{\pi}_0^E$\\
\hline
0.1&0.17805&0.18048&0.23915&0.31011&0.24104&0.26921&0.35413&\textbf{0.14361}&0.15053\\
0.2&0.13587&0.13651&0.18814&0.24789&0.18676&0.21016&0.27621&\textbf{0.11135}&0.11718\\
0.3&0.10685&0.10777&0.15432&0.21188&0.14921&0.16647&0.21339&\textbf{0.09205}&0.09716\\
0.4&0.07843&0.07908&0.11928&0.16674&0.11099&0.12379&0.15134&\textbf{0.07082}&0.07455\\
0.5&0.05107&0.05164&0.08444&0.11871&0.07854&0.08485&0.10310&\textbf{0.04916}&0.05136\\
0.6&\textbf{0.03380}&0.03394&0.05949&0.08309&0.05148&0.05718&0.06382&0.03474&0.03610\\
0.7&0.02079&\textbf{0.01947}&0.03781&0.05357&0.02897&0.03443&0.03338&0.02302&0.02348\\
0.8&0.01239&\textbf{0.00980}&0.01942&0.02862&0.01335&0.01703&0.01446&0.01252&0.01335\\
0.9&0.01015&0.00548&0.00626&0.01812&\textbf{0.00404}&0.00487&0.01527&0.00531&0.00599\\
\hline
\\
\multicolumn{10}{c}{m=500    n=50}\\
\\
\hline
$\pi_0$&$\hat{\pi}_0^B$&$\hat{\pi}_0^L$&$\hat{\pi}_0^A$&$\hat{\pi}_0^P$&$\hat{\pi}_0^H$&$\hat{\pi}_0^D$&$\hat{\pi}_0^S$&$\hat{\pi}_0^U$&$\hat{\pi}_0^E$\\
\hline
0.1&0.10476&0.10384&0.11993&0.11878&0.12536&0.16464&0.25981&\textbf{0.07591}&0.08044\\
0.2&0.08417&0.08283&0.09857&0.10309&0.10171&0.13101&0.20357&\textbf{0.06346}&0.06742\\
0.3&0.06350&0.06302&0.07733&0.08300&0.08049&0.10131&0.15677&\textbf{0.05054}&0.05342\\
0.4&0.04631&0.04689&0.05881&0.06430&0.06192&0.07484&0.11378&\textbf{0.03907}&0.04097\\
0.5&0.03230&0.03253&0.04273&0.05320&0.04358&0.05242&0.07918&\textbf{0.02828}&0.02994\\
0.6&0.01993&0.02028&0.02870&0.03872&0.02923&0.03393&0.05060&\textbf{0.01894}&0.01997\\
0.7&0.01164&\textbf{0.01154}&0.01809&0.02905&0.01717&0.01962&0.02827&0.01169&0.01235\\
0.8&0.00568&\textbf{0.00535}&0.00960&0.01685&0.00800&0.00933&0.01234&0.00616&0.00637\\
0.9&0.00285&\textbf{0.00179}&0.00357&0.00778&0.00211&0.00281&0.00288&0.00227&0.00236\\
\hline
\\
\multicolumn{10}{c}{m=500    n=30}\\
\\
\hline
$\pi_0$&$\hat{\pi}_0^B$&$\hat{\pi}_0^L$&$\hat{\pi}_0^A$&$\hat{\pi}_0^P$&$\hat{\pi}_0^H$&$\hat{\pi}_0^D$&$\hat{\pi}_0^S$&$\hat{\pi}_0^U$&$\hat{\pi}_0^E$\\
\hline
0.1&0.17154&0.16989&0.19307&0.19026&0.19994&0.25450&0.38188&\textbf{0.13731}&0.14537\\
0.2&0.13101&0.13078&0.15195&0.15032&0.16034&0.20034&0.30222&\textbf{0.10898}&0.11499\\
0.3&0.10216&0.10225&0.12106&0.12516&0.12709&0.15509&0.23116&\textbf{0.08773}&0.09260\\
0.4&0.07343&0.07376&0.08950&0.09390&0.09433&0.11371&0.16972&\textbf{0.06571}&0.06922\\
0.5&0.05225&0.05291&0.06610&0.07466&0.06910&0.08056&0.11709&\textbf{0.04910}&0.05143\\
0.6&0.03197&0.03284&0.04306&0.05136&0.04443&0.05143&0.07487&\textbf{0.03163}&0.03339\\
0.7&\textbf{0.01814}&0.01858&0.02646&0.03420&0.02600&0.02952&0.04178&0.01912&0.02013\\
0.8&0.00811&\textbf{0.00798}&0.01303&0.01803&0.01228&0.01353&0.01837&0.00928&0.00964\\
0.9&0.00307&\textbf{0.00224}&0.00426&0.00759&0.00316&0.00380&0.00419&0.00292&0.00311\\
\hline
\\
\multicolumn{10}{c}{m=1000    n=50}\\
\\
\hline
$\pi_0$&$\hat{\pi}_0^B$&$\hat{\pi}_0^L$&$\hat{\pi}_0^A$&$\hat{\pi}_0^P$&$\hat{\pi}_0^H$&$\hat{\pi}_0^D$&$\hat{\pi}_0^S$&$\hat{\pi}_0^U$&$\hat{\pi}_0^E$\\
\hline
0.1&0.10499&0.10350&0.11524&0.10990&0.11971&0.16384&0.26162&\textbf{0.07674}&0.08111\\
0.2&0.08338&0.08261&0.09346&0.09002&0.09743&0.13033&0.20683&\textbf{0.06365}&0.06715\\
0.3&0.06367&0.06285&0.07307&0.07249&0.07601&0.09998&0.15797&\textbf{0.05057}&0.05339\\
0.4&0.04648&0.04610&0.05452&0.05410&0.05752&0.07381&0.11649&\textbf{0.03866}&0.04072\\
0.5&0.03075&0.03067&0.03796&0.04048&0.04022&0.05041&0.07995&\textbf{0.02669}&0.02829\\
0.6&0.01983&0.01998&0.02570&0.02950&0.02697&0.03292&0.05157&\textbf{0.01836}&0.01944\\
0.7&\textbf{0.01095}&0.01122&0.01545&0.01985&0.01629&0.01888&0.02885&0.01123&0.01174\\
0.8&0.00517&\textbf{0.00512}&0.00807&0.01287&0.00785&0.00888&0.01295&0.00570&0.00598\\
0.9&0.00177&\textbf{0.00140}&0.00269&0.00549&0.00199&0.00237&0.00306&0.00181&0.00182\\
\hline
\\
\multicolumn{10}{c}{m=1000   n=30}\\
\\
\hline
$\pi_0$&$\hat{\pi}_0^B$&$\hat{\pi}_0^L$&$\hat{\pi}_0^A$&$\hat{\pi}_0^P$&$\hat{\pi}_0^H$&$\hat{\pi}_0^D$&$\hat{\pi}_0^S$&$\hat{\pi}_0^U$&$\hat{\pi}_0^E$\\
\hline
0.1&0.17211&0.16870&0.18517&0.17670&0.18924&0.25283&0.38573&\textbf{0.13727}&0.14558\\
0.2&0.13502&0.13297&0.14794&0.14215&0.15308&0.19940&0.30311&\textbf{0.11084}&0.11746\\
0.3&0.10272&0.10167&0.11518&0.11092&0.12081&0.15387&0.23327&\textbf{0.08800}&0.09290\\
0.4&0.07445&0.07404&0.08520&0.08305&0.08946&0.11253&0.17028&\textbf{0.06586}&0.06955\\
0.5&0.05065&0.05124&0.06010&0.05882&0.06411&0.07852&0.11867&\textbf{0.04735}&0.04978\\
0.6&0.03181&0.03212&0.03952&0.04192&0.04201&0.05001&0.07559&\textbf{0.03095}&0.03266\\
0.7&\textbf{0.01818}&0.01847&0.02368&0.02732&0.02457&0.02877&0.04257&0.01883&0.01982\\
0.8&\textbf{0.00764}&0.00790&0.01130&0.01430&0.01153&0.01288&0.01862&0.00883&0.00923\\
0.9&0.00230&\textbf{0.00209}&0.00369&0.00561&0.00308&0.00356&0.00454&0.00272&0.00284\\
\hline
\end{longtable}
\begin{longtable}{cccccccccc}
\caption{Bias of estimators under different simulation settings}\\
\endfirsthead
\hline
\\
\multicolumn{10}{c}{m=100     n=50}\\
\\
\hline
$\pi_0$&$\hat{\pi}_0^B$&$\hat{\pi}_0^L$&$\hat{\pi}_0^A$&$\hat{\pi}_0^P$&$\hat{\pi}_0^H$&$\hat{\pi}_0^D$&$\hat{\pi}_0^S$&$\hat{\pi}_0^U$&$\hat{\pi}_0^E$\\
\hline
0.1&0.31797&0.32359&0.38310&0.42538&0.38330&0.41589&0.47213&\textbf{0.26626}&0.27415\\
0.2&0.27582&0.28167&0.34204&0.37872&0.34134&0.36922&0.41578&\textbf{0.23668}&0.24322\\
0.3&0.24167&0.24807&0.31109&0.36059&0.30433&0.32781&0.36001&\textbf{0.21368}&0.21929\\
0.4&0.19575&0.20278&0.26797&0.30633&0.25931&0.28144&0.30659&\textbf{0.18129}&0.18507\\
0.5&0.15647&0.16605&0.22937&0.26073&0.21551&0.23678&0.25026&\textbf{0.15205}&0.15557\\
0.6&\textbf{0.11245}&0.12446&0.19048&0.20392&0.17702&0.19167&0.19396&0.12210&0.12279\\
0.7&\textbf{0.07002}&0.08432&0.14652&0.13948&0.12866&0.14662&0.14241&0.09114&0.09049\\
0.8&\textbf{0.02780}&0.04031&0.09984&0.08517&0.07979&0.09785&0.07951&0.05487&0.05552\\
0.9&-0.01468&-0.00145&0.04520&0.01675&0.03176&0.04790&\textbf{0.00685}&0.01467&0.01804\\
\hline
\\
\multicolumn{10}{c}{m=100     n=30}\\
\\
\hline
$\pi_0$&$\hat{\pi}_0^B$&$\hat{\pi}_0^L$&$\hat{\pi}_0^A$&$\hat{\pi}_0^P$&$\hat{\pi}_0^H$&$\hat{\pi}_0^D$&$\hat{\pi}_0^S$&$\hat{\pi}_0^U$&$\hat{\pi}_0^E$\\
\hline
0.1&0.40921&0.41469&0.48010&0.51813&0.48405&0.51552&0.59179&\textbf{0.36447}&0.37342\\
0.2&0.35312&0.35733&0.42276&0.45682&0.42419&0.45450&0.52228&\textbf{0.31610}&0.32490\\
0.3&0.30829&0.31501&0.38072&0.41417&0.37860&0.40367&0.45877&\textbf{0.28516}&0.29311\\
0.4&0.25974&0.26680&0.33277&0.36115&0.32518&0.34703&0.38535&\textbf{0.24740}&0.25444\\
0.5&0.19958&0.20818&0.27419&0.28783&0.27067&0.28536&0.31701&\textbf{0.19888}&0.20308\\
0.6&\textbf{0.14907}&0.16170&0.22607&0.22758&0.21817&0.23202&0.24751&0.16059&0.16341\\
0.7&\textbf{0.09981}&0.11304&0.17513&0.16738&0.16008&0.17627&0.17620&0.12031&0.12218\\
0.8&\textbf{0.04910}&0.06175&0.11725&0.09819&0.10353&0.11918&0.10389&0.07545&0.07787\\
0.9&-0.00710&0.00558&0.04875&0.01317&0.04539&0.05427&\textbf{0.00016}&0.02355&0.02517\\
\hline
\\
\multicolumn{10}{c}{m=500     n=50}\\
\\
\hline
$\pi_0$&$\hat{\pi}_0^B$&$\hat{\pi}_0^L$&$\hat{\pi}_0^A$&$\hat{\pi}_0^P$&$\hat{\pi}_0^H$&$\hat{\pi}_0^D$&$\hat{\pi}_0^S$&$\hat{\pi}_0^U$&$\hat{\pi}_0^E$\\
\hline
0.1&0.32001&0.31906&0.34319&0.33376&0.35100&0.40495&0.50912&\textbf{0.27180}&0.27978\\
0.2&0.28597&0.28431&0.31041&0.30792&0.31572&0.36097&0.45057&\textbf{0.24778}&0.25546\\
0.3&0.24691&0.24686&0.27423&0.27131&0.28046&0.31729&0.39537&\textbf{0.22032}&0.22666\\
0.4&0.20876&0.21182&0.23811&0.23289&0.24540&0.27246&0.33674&\textbf{0.19300}&0.19760\\
0.5&0.17210&0.17484&0.20171&0.20598&0.20501&0.22760&0.28081&\textbf{0.16278}&0.16762\\
0.6&\textbf{0.13073}&0.13535&0.16294&0.16440&0.16703&0.18241&0.22422&0.13122&0.13471\\
0.7&\textbf{0.09405}&0.09796&0.12595&0.13288&0.12582&0.13771&0.16730&0.10035&0.10317\\
0.8&\textbf{0.05306}&0.06105&0.08732&0.08243&0.08495&0.09320&0.11019&0.06881&0.06984\\
0.9&\textbf{0.01188}&0.01983&0.04440&0.02821&0.03907&0.04665&0.05224&0.03225&0.03345\\
\hline
\\
\multicolumn{10}{c}{m=500     n=30}\\
\\
\hline
$\pi_0$&$\hat{\pi}_0^B$&$\hat{\pi}_0^L$&$\hat{\pi}_0^A$&$\hat{\pi}_0^P$&$\hat{\pi}_0^H$&$\hat{\pi}_0^D$&$\hat{\pi}_0^S$&$\hat{\pi}_0^U$&$\hat{\pi}_0^E$\\
\hline
0.1&0.41076&0.40914&0.43653&0.42595&0.44432&0.50380&0.61756&\textbf{0.36735}&0.37804\\
0.2&0.35766&0.35811&0.38639&0.37454&0.39721&0.44681&0.54932&\textbf{0.32665}&0.33546\\
0.3&0.31146&0.31583&0.34431&0.33814&0.35327&0.39292&0.48035&\textbf{0.29247}&0.30045\\
0.4&0.26532&0.26728&0.29518&0.28790&0.30377&0.33625&0.41152&\textbf{0.25238}&0.25908\\
0.5&0.22182&0.22530&0.25288&0.25265&0.25981&0.28264&0.34170&\textbf{0.21720}&0.22236\\
0.6&\textbf{0.16941}&0.17494&0.20165&0.19922&0.20705&0.22532&0.27311&0.17273&0.17738\\
0.7&\textbf{0.12334}&0.12903&0.15581&0.15243&0.15764&0.16995&0.20384&0.13240&0.13601\\
0.8&\textbf{0.06921}&0.07810&0.10479&0.09262&0.10642&0.11344&0.13487&0.08801&0.08989\\
0.9&\textbf{0.02120}&0.03049&0.05320&0.03381&0.05072&0.05692&0.06372&0.04323&0.04451\\
\hline
\\
\multicolumn{10}{c}{m=1000    n=50}\\
\\
\hline
$\pi_0$&$\hat{\pi}_0^B$&$\hat{\pi}_0^L$&$\hat{\pi}_0^A$&$\hat{\pi}_0^P$&$\hat{\pi}_0^H$&$\hat{\pi}_0^D$&$\hat{\pi}_0^S$&$\hat{\pi}_0^U$&$\hat{\pi}_0^E$\\
\hline
0.1&0.32191&0.31990&0.33754&0.32531&0.34401&0.40438&0.51120&\textbf{0.27514}&0.28279\\
0.2&0.28618&0.28538&0.30360&0.29287&0.31004&0.36057&0.45452&\textbf{0.25022}&0.25697\\
0.3&0.24920&0.24814&0.26785&0.25984&0.27357&0.31567&0.39715&\textbf{0.22255}&0.22861\\
0.4&0.21178&0.21148&0.23051&0.22102&0.23740&0.27099&0.34099&\textbf{0.19385}&0.19892\\
0.5&0.17081&0.17151&0.19152&0.18628&0.19801&0.22376&0.28245&\textbf{0.16045}&0.16513\\
0.6&0.13473&0.13672&0.15628&0.15190&0.16138&0.18049&0.22677&\textbf{0.13214}&0.13588\\
0.7&\textbf{0.09627}&0.10067&0.11938&0.11428&0.12498&0.13620&0.16951&0.10225&0.10433\\
0.8&\textbf{0.05913}&0.06474&0.08372&0.08100&0.08540&0.09259&0.11342&0.07074&0.07250\\
0.9&\textbf{0.01746}&0.02453&0.04126&0.02885&0.04016&0.04526&0.05464&0.03415&0.03460\\
\hline
\\
\multicolumn{10}{c}{m=1000    n=30}\\
\\
\hline
$\pi_0$&$\hat{\pi}_0^B$&$\hat{\pi}_0^L$&$\hat{\pi}_0^A$&$\hat{\pi}_0^P$&$\hat{\pi}_0^H$&$\hat{\pi}_0^D$&$\hat{\pi}_0^S$&$\hat{\pi}_0^U$&$\hat{\pi}_0^E$\\
\hline
0.1&0.41304&0.40911&0.42867&0.41539&0.43326&0.50248&0.62086&\textbf{0.36896}&0.37996\\
0.2&0.36520&0.36264&0.38275&0.37060&0.38935&0.44614&0.55034&\textbf{0.33119}&0.34091\\
0.3&0.31785&0.31659&0.33720&0.32467&0.34546&0.39183&0.48277&\textbf{0.29476}&0.30286\\
0.4&0.26955&0.26947&0.28945&0.27789&0.29695&0.33496&0.41242&\textbf{0.25458}&0.26157\\
0.5&0.22105&0.22333&0.24234&0.22945&0.25082&0.27961&0.34426&\textbf{0.21527}&0.22073\\
0.6&\textbf{0.17349}&0.17571&0.19540&0.18785&0.20252&0.22290&0.27469&0.17335&0.17807\\
0.7&\textbf{0.12834}&0.13140&0.14966&0.14535&0.15364&0.16857&0.20605&0.13375&0.13746\\
0.8&\textbf{0.07671}&0.08265&0.10045&0.09173&0.10436&0.11199&0.13612&0.08991&0.09194\\
0.9&\textbf{0.02926}&0.03555&0.05228&0.03904&0.05198&0.05681&0.06686&0.04585&0.04702\\
\hline
\end{longtable}

\begin{longtable}{cccccccccc}
\caption{Power of adaptive algorithm with the mentioned estimators under different simulation settings}\\
\endfirsthead
\hline
\\
\multicolumn{10}{c}{m=100    n=50}\\
\\
\hline
$\pi_0$&$\hat{\pi}_0^B$&$\hat{\pi}_0^L$&$\hat{\pi}_0^A$&$\hat{\pi}_0^P$&$\hat{\pi}_0^H$&$\hat{\pi}_0^D$&$\hat{\pi}_0^S$&$\hat{\pi}_0^U$&$\hat{\pi}_0^E$\\
\hline
0.1&0.51576&0.51366&0.50023&0.49972&0.50061&0.49263&0.48361&\textbf{0.53268}&0.52959\\
0.2&0.45198&0.45007&0.44014&0.44012&0.44032&0.43484&0.42784&\textbf{0.46245}&0.45998\\
0.3&0.38891&0.38734&0.37925&0.37867&0.37947&0.37612&0.37203&\textbf{0.39413}&0.39321\\
0.4&0.33100&0.32962&0.32326&0.32348&0.32364&0.32142&0.31871&\textbf{0.33294}&0.33221\\
0.5&0.27511&0.27394&0.26927&0.26938&0.27028&0.26849&0.26709&\textbf{0.27537}&0.27520\\
0.6&\textbf{0.22298}&0.22190&0.21868&0.21987&0.21914&0.21841&0.21823&0.22213&0.22224\\
0.7&\textbf{0.17235}&0.17146&0.16950&0.17052&0.17008&0.16928&0.16960&0.17152&0.17150\\
0.8&\textbf{0.12922}&0.12893&0.12766&0.12832&0.12804&0.12768&0.12858&0.12876&0.12863\\
0.9&\textbf{0.08675}&0.08665&0.08620&0.08657&0.08635&0.08619&0.08639&0.08648&0.08649\\
\hline
\\
\multicolumn{10}{c}{m=100    n=30}\\
\\
\hline
$\pi_0$&$\hat{\pi}_0^B$&$\hat{\pi}_0^L$&$\hat{\pi}_0^A$&$\hat{\pi}_0^P$&$\hat{\pi}_0^H$&$\hat{\pi}_0^D$&$\hat{\pi}_0^S$&$\hat{\pi}_0^U$&$\hat{\pi}_0^E$\\
\hline
0.1&0.38419&0.38218&0.37046&0.37226&0.36902&0.36328&0.35211&\textbf{0.39539}&0.39322\\
0.2&0.33595&0.33430&0.32482&0.32538&0.32424&0.31956&0.31067&\textbf{0.34434}&0.34190\\
0.3&0.28416&0.28247&0.27525&0.27654&0.27484&0.27243&0.26731&\textbf{0.28735}&0.28626\\
0.4&0.24422&0.24312&0.23756&0.23808&0.23800&0.23617&0.23354&\textbf{0.24566}&0.24502\\
0.5&\textbf{0.20542}&0.20445&0.20020&0.20224&0.20040&0.19922&0.19788&0.20522&0.20514\\
0.6&\textbf{0.16578}&0.16491&0.16236&0.16355&0.16250&0.16202&0.16150&0.16526&0.16515\\
0.7&\textbf{0.13367}&0.13304&0.13114&0.13206&0.13166&0.13119&0.13102&0.13280&0.13287\\
0.8&\textbf{0.10108}&0.10070&0.09975&0.10043&0.10004&0.09971&0.09997&0.10051&0.10040\\
0.9&0.07416&\textbf{0.07478}&0.07446&0.07472&0.07455&0.07446&0.07476&0.07462&0.07464\\
\hline
\\
\multicolumn{10}{c}{m=500    n=50}\\
\\
\hline
$\pi_0$&$\hat{\pi}_0^B$&$\hat{\pi}_0^L$&$\hat{\pi}_0^A$&$\hat{\pi}_0^P$&$\hat{\pi}_0^H$&$\hat{\pi}_0^D$&$\hat{\pi}_0^S$&$\hat{\pi}_0^U$&$\hat{\pi}_0^E$\\
\hline
0.1&0.51141&0.51148&0.50638&0.50987&0.50488&0.49407&0.47863&\textbf{0.52362}&0.52133\\
0.2&0.44682&0.44691&0.44268&0.44474&0.44179&0.43485&0.42408&\textbf{0.45404}&0.45250\\
0.3&0.38463&0.38449&0.38098&0.38269&0.38012&0.37581&0.36819&\textbf{0.38837}&0.38746\\
0.4&0.23832&0.32785&0.32548&0.32686&0.32477&0.32237&0.31744&\textbf{0.32981}&0.32933\\
0.5&0.27234&0.27212&0.27036&0.27071&0.27018&0.26856&0.26553&\textbf{0.27299}&0.27262\\
0.6&\textbf{0.22074}&0.22042&0.21910&0.21954&0.21889&0.21820&0.21639&0.22063&0.22046\\
0.7&\textbf{0.17164}&0.17147&0.17068&0.17079&0.17065&0.17029&0.16945&0.17137&0.17131\\
0.8&\textbf{0.12626}&0.12608&0.12563&0.12586&0.12568&0.12550&0.12523&0.12591&0.12592\\
0.9&\textbf{0.08521}&0.08514&0.08496&0.08509&0.08502&0.08495&0.08491&0.08505&0.08505\\
\hline
\\
\multicolumn{10}{c}{m=500    n=30}\\
\\
\hline
$\pi_0$&$\hat{\pi}_0^B$&$\hat{\pi}_0^L$&$\hat{\pi}_0^A$&$\hat{\pi}_0^P$&$\hat{\pi}_0^H$&$\hat{\pi}_0^D$&$\hat{\pi}_0^S$&$\hat{\pi}_0^U$&$\hat{\pi}_0^E$\\
\hline
0.1&0.37771&0.37791&0.37289&0.37635&0.37155&0.36157&0.34710&\textbf{0.38681}&0.38444\\
0.2&0.32869&0.32852&0.32429&0.32758&0.32299&0.31638&0.30550&\textbf{0.33358}&0.33215\\
0.3&0.28103&0.28084&0.27774&0.27968&0.27663&0.27267&0.26538&\textbf{0.28363}&0.28261\\
0.4&0.23836&0.23804&0.23581&0.23730&0.23499&0.23241&0.22761&\textbf{0.23948}&0.23887\\
0.5&0.19846&0.19814&0.19657&0.19711&0.19609&0.19490&0.19194&\textbf{0.19871}&0.19837\\
0.6&\textbf{0.16168}&0.16136&0.16032&0.16090&0.16004&0.15940&0.15781&0.16146&0.16130\\
0.7&\textbf{0.12848}&0.12834&0.12776&0.12814&0.12765&0.12739&0.12673&0.12823&0.12821\\
0.8&\textbf{0.09789}&0.09768&0.09733&0.09762&0.09729&0.09719&0.09697&0.09752&0.09751\\
0.9&\textbf{0.07208}&0.07202&0.07191&0.07205&0.07192&0.07191&0.07118&0.07196&0.07196\\
\hline
\\
\multicolumn{10}{c}{m=1000    n=50}\\
\\
\hline
$\pi_0$&$\hat{\pi}_0^B$&$\hat{\pi}_0^L$&$\hat{\pi}_0^A$&$\hat{\pi}_0^P$&$\hat{\pi}_0^H$&$\hat{\pi}_0^D$&$\hat{\pi}_0^S$&$\hat{\pi}_0^U$&$\hat{\pi}_0^E$\\
\hline
0.1&0.51073&0.51112&0.50739&0.51117&0.50607&0.49463&0.47904&\textbf{0.52202}&0.52001\\
0.2&0.44598&0.44602&0.44312&0.44551&0.44209&0.43451&0.42316&\textbf{0.45243}&0.45114\\
0.3&0.38507&0.38512&0.38272&0.38429&0.38192&0.37702&0.36910&\textbf{0.38856}&0.38772\\
0.4&0.32640&0.32639&0.32461&0.32593&0.32391&0.32096&0.31571&\textbf{0.32820}&0.32766\\
0.5&0.27346&0.27333&0.27200&0.27267&0.27149&0.26984&0.26646&\textbf{0.27411}&0.27379\\
0.6&0.21987&0.21981&0.21889&0.21932&0.21863&0.21779&0.21583&\textbf{0.22000}&0.21982\\
0.7&\textbf{0.17145}&0.17127&0.17071&0.17104&0.17054&0.17019&0.16926&0.17121&0.17115\\
0.8&\textbf{0.12537}&0.12527&0.12498&0.12513&0.12491&0.12482&0.12451&0.12515&0.12515\\
0.9&\textbf{0.08480}&0.08476&0.08465&0.08476&0.08465&0.08462&0.08455&0.08468&0.08467\\
\hline
\\
\multicolumn{10}{c}{m=1000     n=30}\\
\\
\hline
$\pi_0$&$\hat{\pi}_0^B$&$\hat{\pi}_0^L$&$\hat{\pi}_0^A$&$\hat{\pi}_0^P$&$\hat{\pi}_0^H$&$\hat{\pi}_0^D$&$\hat{\pi}_0^S$&$\hat{\pi}_0^U$&$\hat{\pi}_0^E$\\
\hline
0.1&0.37585&0.37656&0.37298&0.37621&0.37225&0.36103&0.34592&\textbf{0.38459}&0.38227\\
0.2&0.32820&0.32848&0.32553&0.32791&0.32459&0.31710&0.30628&\textbf{0.33340}&0.33182\\
0.3&0.28057&0.28062&0.27842&0.28041&0.27752&0.27283&0.26529&\textbf{0.28319}&0.28229\\
0.4&0.23860&0.23858&0.23692&0.23828&0.23633&0.23335&0.22818&\textbf{0.23978}&0.23919\\
0.5&0.19823&0.19802&0.19692&0.19787&0.19648&0.19492&0.19173&\textbf{0.19849}&0.19819\\
0.6&\textbf{0.16158}&0.16144&0.16074&0.16128&0.16041&0.15972&0.15798&0.16153&0.16135\\
0.7&\textbf{0.12728}&0.12720&0.12674&0.12706&0.12666&0.12634&0.12556&0.12714&0.12705\\
0.8&\textbf{0.09790}&0.09783&0.09761&0.09776&0.09754&0.09744&0.09713&0.09773&0.09770\\
0.9&\textbf{0.07152}&0.07148&0.07142&0.07149&0.07143&0.07139&0.07136&0.07144&0.07144\\
\hline
\end{longtable}
\end{small}
\begin{figure}[h!]
\begin{center}
\includegraphics[height=6in,width=6in,angle=0]{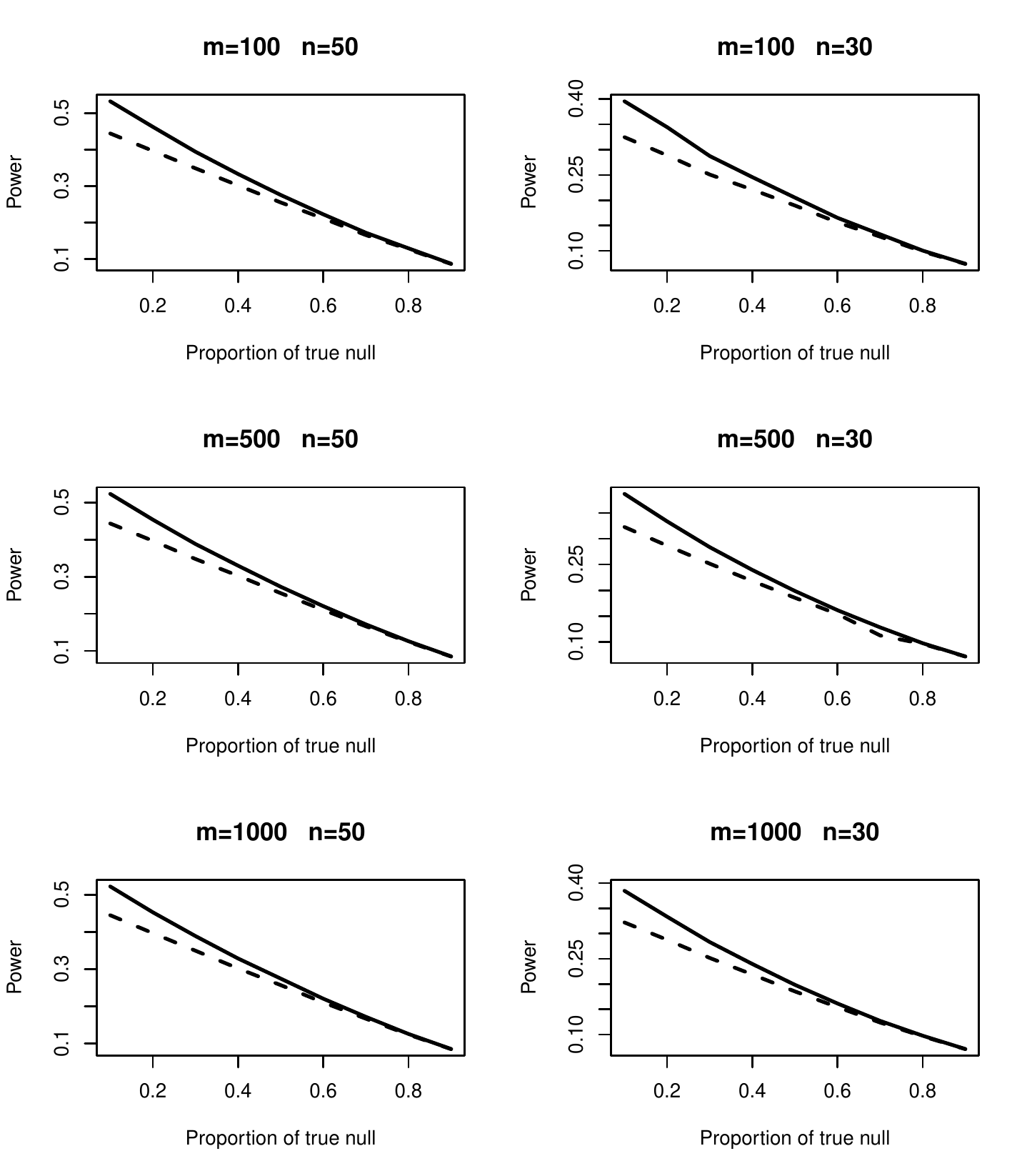}
\caption{Plot showing power of non-adaptive BH algorithm (dashed line) and adaptive BH algorith with $\hat{\pi}_0^U$ (solid line) as functions of $\pi_0$ under different simulation settings.}
\end{center}
\end{figure}
\end{singlespacing}

\section{Data analysis}
For the case-study we have considered the real life synthetic data set used by Gupta et al. (2014, 2017) in connection with reliability and warranty studies. The detailed description of the data is available there and we report only the relevant part of it, which is required in the present study. The date of failure of a particular component of  automobiles along with the mileage at failure as reflected through the odometer readings are available. Although the entire data set cover two disjoint geographical regions, as reported in Gupta et al. (2017), they may be further subdivided into failures corresponding to seven sub-regions, termed as zones. Owing to data confidentiality issue, let us number them from 1 to 7. We have considered failure data corresponding to a particular year and the mileage figures of the failures in successive calendar months across the zones as the response variable. The twelve calendar months are recorded as JAN, FEB, MAR, APR, MAY, JUN, JUL, AUG, SEP, OCT, NOV and DEC. Thus the entire dataset is inherently segmented by $7$ different zones and $12$ different months in a year. In other words, our synthetic dataset contains mileage at failure for $84$ different month$\times$zone segments. In total, little less than 3000 component failures have been reported in the year considered with varying number of warranty claims over the month$\times$zone segments. Thus on an average around 35 failures are reported in each segment. In line with the discussion in section 1, here we are primarily interested in identifying the segments which have significantly better or poor performance in terms of mileage at failure in comparison to a bench mark. Thus, appropriate hypotheses are needed to be formed  and tested separately for each of the segments making way for application of adaptive FDR-controlling algorithms.\\

To validate our model-assumption we perform Kolmogorov-Smirnov's test with empirical $p$-value for exponential distribution (using R-function \begin{slshape}ks.exp.test\end{slshape} available in \begin{slshape}exptest\end{slshape} package) and find that, at level $0.05$ exponentiality fails for only $18$ out of $84$ segments whereas at level $0.01$ only $7$ rejections are there. For visual display of fit, we present QQ-plot of some randomly chosen segments in Figure 2. As the sample sizes for most of the segments are moderate, we also check normality applying Shapiro-Wilks' test (using default R-function \begin{slshape}shapiro.test\end{slshape}). At level $0.05$, $59$ out of $84$ hypotheses gets rejected and at level $0.01$, the number is $42$. The first line of information justifies applicability of the model based estimators for $\pi_0$ discussed in this article whereas the results from the normality test demonstrate the necessity of the modifications achieved through this work over the existing related works, usually done under normal model.   \\   

\begin{figure}[h!]
\begin{center}
\includegraphics[height=6in,width=6in,angle=0]{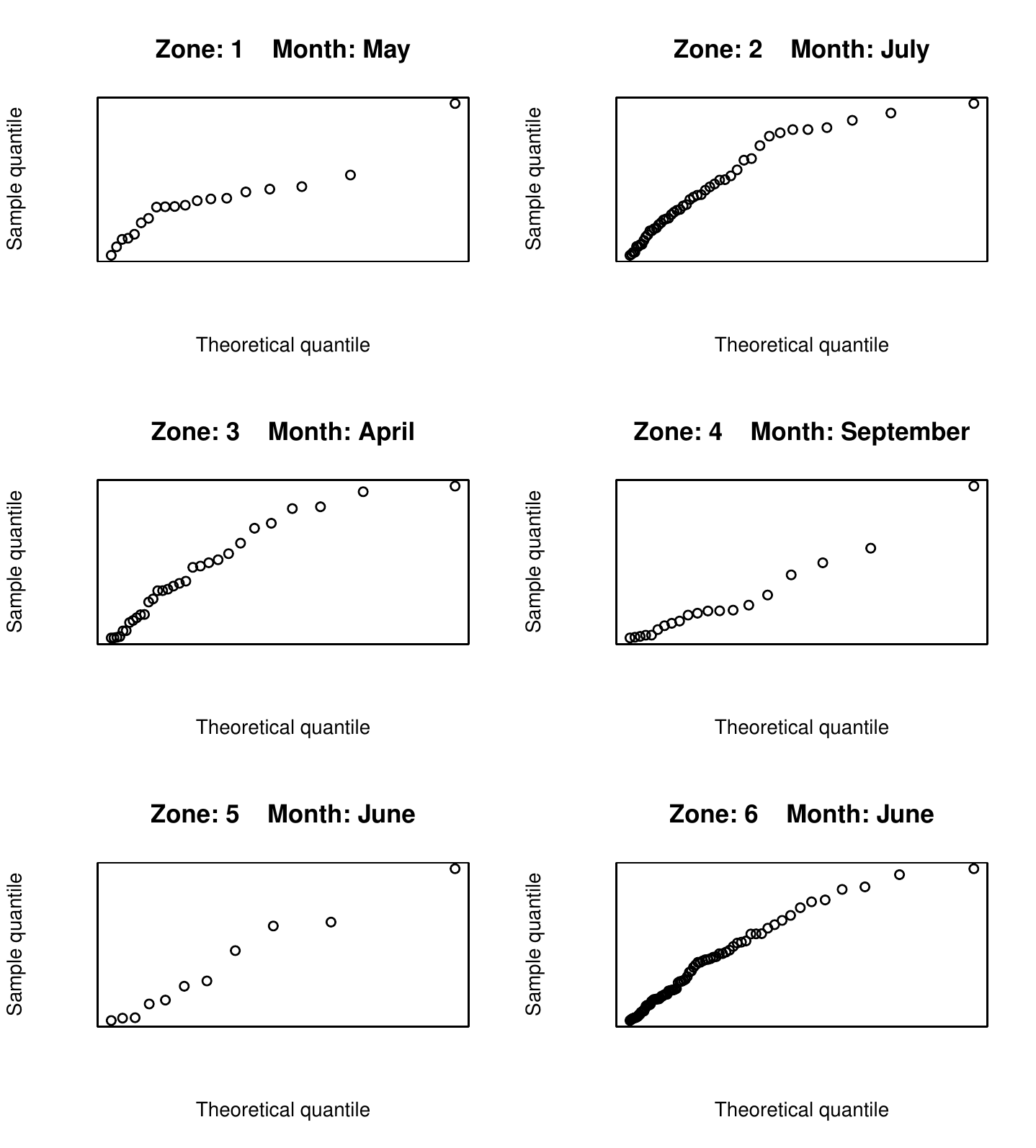}
\caption{QQplot for some randomly chosen segments from the synthetic data, showcasing the exponential fitting. }
\end{center}
\end{figure}

Now we consider framing of the appropriate hypotheses. We assume that the mileages at failure for the $i$-th segment to be exponentially distributed with mean $\theta_i$ miles. Thus $\theta_i$'s are the mean mileage to failure (MMTF) for the $i$-th. segment, a quantity similar to mean time to failure (MTTF) in terms of the response variable `mileage', for $i = 1, 2, ..., 84$. We consider, as an indicator of the bench mark, the MMTF of the entire dataset as our null hypothesis point, approximately given by $\theta_0=10973$ miles. This value as an benchmark seems to be justified, as the warranty mileage limit for the data base is $36000$ miles and it is well known that such failure data are usually positively skewed. According to the research question we then simultaneously test the following hypotheses:
\begin{equation}\label{eq:datahyp}
H_{0i}:\theta_i=\theta_0 \quad\quad\textrm{versus}\quad\quad H_{1i}:\theta_i\neq\theta_0\quad\textrm{for}\quad i=1, 2, ..., 84.  
\end{equation}
The two-sided choice of the alternative hypotheses at all the segments needs clarification. In the absence of any prior knowledge about the functioning of the component, it is not possible to mark any segment to be better/worse than the overall benchmark in terms of MMTF. As a result, to be on the safe side, we have suggested the alternative hypotheses at all segments to be two-sided. This is very common in multiple testing situations. As an example, in microarray data analysis, the samples used as a reference are called control samples. The other samples exhibiting different phenotypic status are called treated samples. The gene expression levels among these groups may be different. To identify whether a gene is differentially expressed or not, we fix two-sided alternative (Chen et al. 2007).\\

Likelihood ratio tests are performed for each of the hypotheses after scaling the original observations by $\theta_0$, maintaining equivalence of the test and the corresponding $p$-values along with effect sizes for each test are stored for further use. A file \begin{slshape}auto\_details.txt\end{slshape} available as supplementary material provides details under the following heads:
\begin{itemize}
\item \textit{segment}: This column provides serial number of the segment, $1$ to $84$ such that segment $i$ is for $i$-th month of zone 1 for $i=1(1)12$, segment $12+i$ is for $i$-th month of zone 2 for $i=1(1)12$ and so on for the $7$ zones in order as mentioned in the first paragraph of this section.
\item \textit{n}: Provides available sample size for each segment.
\item \textit{pval}: Provides the obtained $p$-value corresponding to common likelihood ratio test performed for each segment.
\item \textit{del}: Provides maximum likelihood estimates of the effect sizes corresponding to each test.  
\end{itemize}
These array of values can be readily fed into Algorithms 1, 2 and 3 to get $\hat{\pi}_0^U$, $\hat{\pi}_0^E$ and the list of rejected hypotheses when adaptive FDR-controlling algorithm is applied with different $\pi_0$-estimates. Estimate of $\pi_0$ along with the corresponding list of rejected hypotheses using the estimators already mentioned in this article are also reported. The estimated $\pi_0$-values using different estimators are reported in Table 4.\\ 
\begin{singlespacing}
\begin{small}
\begin{table}[h!]\label{tab:datapi0estim}
\caption{Estimate of $\pi_0$ using different estimators under exponential model assumption for the synthetic data.}
\center
\begin{tabular}{cccccc}
\hline
\\
$\hat{\pi}_0^B$ & 0.5555 & $\hat{\pi}_0^L$ & 0.5565 & $\hat{\pi}_0^A$ & 0.5761\\
\\
$\hat{\pi}_0^P$ & 0.6074 & $\hat{\pi}_0^H$ & 0.5555 & $\hat{\pi}_0^D$ & 0.6662\\
\\
$\hat{\pi}_0^S$ & 0.8065 & $\hat{\pi}_0^U$ & 0.4842 & $\hat{\pi}_0^E$ & 0.5096\\
\\
\hline
\end{tabular}
\end{table} 
\end{small}
\end{singlespacing}

In Table 5 we indicate the segments that are found to be significantly different from the average in terms of the mean mileage at of the designated component failure when adaptive BH-algorithm with different $\pi_0$-estimators and non-adaptive(N/A) BH-algorithm are applied to control FDR at level $q=0.05,0.1$. For visual display, we plot the adjusted $p$-values for non-adaptive Benjamini-Hochberg algorithm and its adaptive version using $\hat{\pi}_0^U$ with cut-off $q=0.1$ in Figure 3. From Figure 3 and Table 5, it is evident that, the adaptive BH-algorithm using the proposed methods has the ability to identify a larger number of segments with significant variation from benchmark by controlling the FDR at the same level, compared to the non-adaptive BH-algorithm as well as adaptive BH-algorithm using existent estimators for $\pi_0$. \\

\begin{singlespacing}
\begin{table}[h!]\label{tab:rejindicator}
\caption{Significantly different segments identified by adaptive-BH algorithm with different estimators for $\pi_0$ for the synthetic data.}
\resizebox{!}{0.325\textwidth}
{
\begin{tabular}{ccccccccccccccc}
\\
\hline
\\
Segment & Zone & Month & mean:$\bar{X}$ & 95\%CI($\theta$)  & $\hat{\pi}_0^B$ & $\hat{\pi}_0^L$ & $\hat{\pi}_0^A$ & $\hat{\pi}_0^P$ & $\hat{\pi}_0^H$ & $\hat{\pi}_0^D$ & $\hat{\pi}_0^S$ & $\hat{\pi}_0^U$ & $\hat{\pi}_0^E$ & N/A\\
\\
\hline
\\
16&2&APR&1.44&(1.12,1.93)&2&2&2&2&2&2&2&2&2&2\\
18&2&JUN&1.29&(1.01,1.70)&0&0&0&0&0&1&0&1&0&0\\
20&2&AUG&1.47&(1.11,2.04)&2&2&2&2&2&2&1&2&2&0\\
21&2&SEP&1.33&(1.01,1.82)&0&0&0&0&0&0&0&1&0&0\\
23&2&NOV&1.69&(1.30,2.29)&2&2&2&2&2&2&2&2&2&2\\
38&4&MAR&0.55&(0.37,0.92)&1&1&1&1&1&0&0&1&1&0\\
39&4&APR&0.50&(0.34,0.79)&2&2&2&2&2&2&2&2&2&2\\
40&4&MAY&0.58&(0.40,0.92)&1&1&1&1&1&0&0&1&1&0\\
41&4&JUN&0.63&(0.44,0.95)&0&0&0&0&0&0&0&1&1&0\\
46&4&OCT&0.64&(0.45,0.98)&0&0&0&0&0&0&0&1&0&0\\
52&5&APR&1.77&(1.10,3.33)&1&1&1&1&1&1&0&1&1&0\\
56&5&AUG&1.82&(1.14,3.32)&1&1&1&1&1&1&1&2&2&1\\
62&6&FEB&0.66&(0.53,0.84)&2&2&2&2&2&2&2&2&2&2\\
63&6&MAR&0.60&(0.48,0.78)&2&2&2&2&2&2&2&2&2&2\\
64&6&APR&1.27&(1.03,1.61)&1&1&1&1&1&0&0&1&1&0\\
69&6&SEP&0.67&(0.54,0.85)&2&2&2&2&2&2&2&2&2&2\\
74&7&FEB&0.20&(0.11,0.51)&2&2&2&2&2&2&2&2&2&2\\
76&7&APR&0.48&(0.20,0.48)&0&0&0&0&0&0&0&1&0&0\\
\\
\hline
\end{tabular}}
\footnotemark{Values in fourth and fifth columns are reported after scaling the original variable by $\theta_0=10973.30$. For columns $6$ to $15$, $2$ indicates that the segments are found significant for $q=0.05$ and trivially for $q=0.1$ both while $1$ indicates same only for $q=0.1$ and $0$ indicates negation of the previous two statements. The segments not reported in this table are not found to be significantly different from the the overall average, taken as the null hypothesis point.}
\end{table}

\begin{figure}[h!]
\begin{center}
\includegraphics[height=4.5in,width=6in,angle=0]{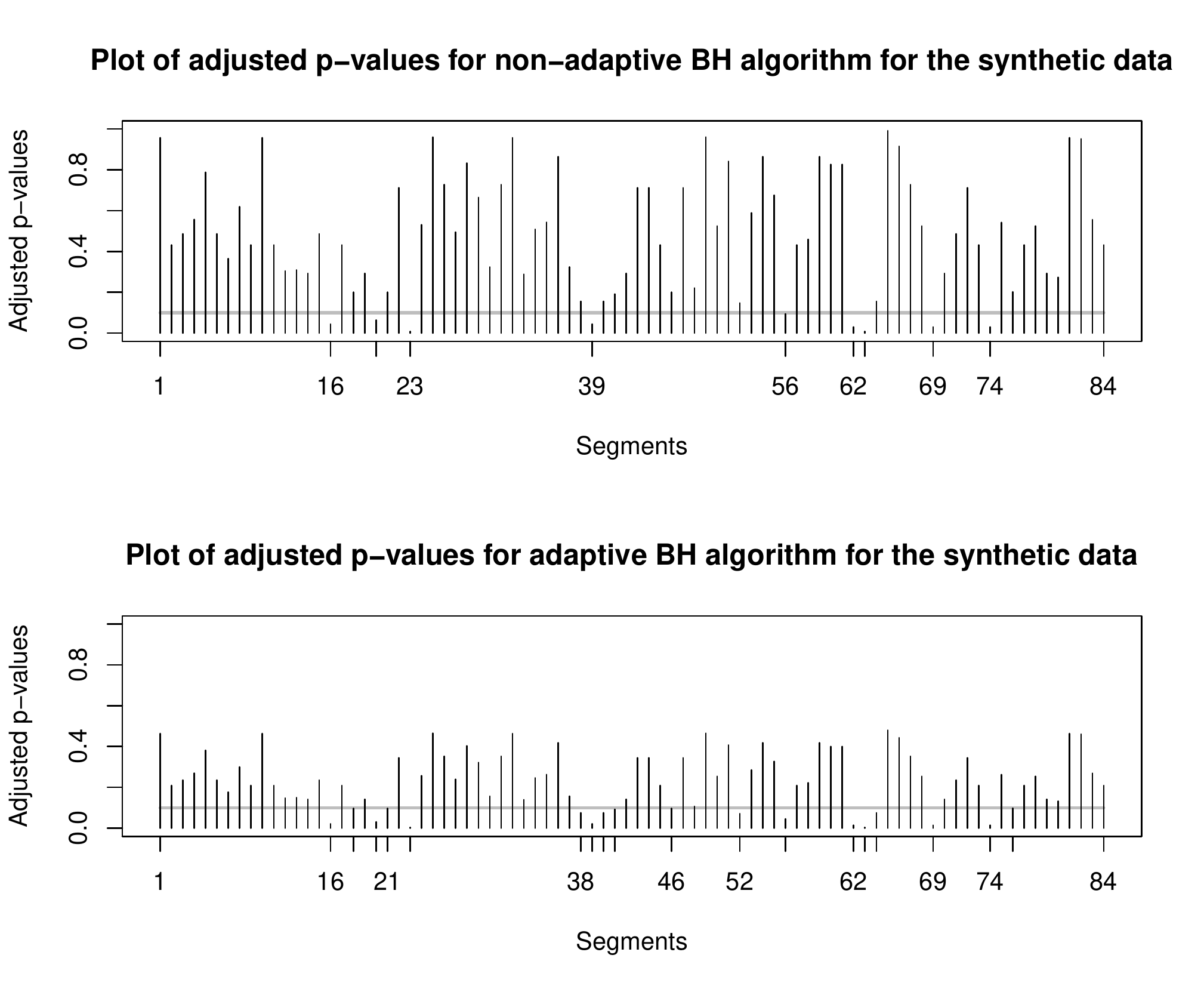}
\caption{Comparative plot of adjusted p-values for the synthetic data with $q=0.1$ as cut-off for visualising gain in power for adaptive BH algorithm with $\hat{\pi}_0^U$. }
\end{center}
\end{figure}
\end{singlespacing}
From the domain knowledge (not to be mentioned explicitly, owing to confidentiality issue), it is known that the functioning of the automobile component under consideration is likely to be influenced by the climate condition, reflected through the effect of the month, as well as by the effect of the zone of their usual operation. The effect of climate on the functioning of the automobiles is well known and has also been reported in Lawless (1998). For simplicity and demonstration purpose, we assume that each automobile is used only in the designated zone where the failure is reported. Although, we have used the two sided alternative, as being done in any multiple testing problem, the corresponding confidence interval falling entirely below or above of the scaled null hypothesis point of `unity', indicates the actual one-sided alternative for which the respective significance appears. Thus the MMTFs of zone 4 are consistently and significantly smaller than the benchmark value (the null hypothesis point) indicating usage related adverse problem of the automobiles and this problem is persistent in the first or second quarter of the year indicating a transition from colder to warmer climate or the fourth quarter of the year indicating the transition from warmer to colder climate . Interestingly for zone 2, exactly the converse situation is prevailed and this seemingly high MMTF might not be due to the climate condition and on the contrary may be attributed to better usage scenario. For zone 5, better usage scenario is evident at least in two months, although weather related issues might not be associated with such improvement. The findings for zone 6 is heavily dependent on climate condition especially during the advent of spring where a significant decrease in MMTF is identified followed by significant increase in MMTF just after. Again during the fall a significant decrease in MMTF is found establishing the climate dependence of failure data. For zone 7, climate plays an adverse role during the end of the winter and start of the summer. The data corresponding to remaining two zones, do not reveal any deviation from the usual usage pattern and/or are not affected by extreme climate conditions. It is to be noted that for almost all the zones, the month of April becomes significant concerning either betterment or worsening of the scenario in comparison to the benchmark. On the other hand the two months viz. December and January never become markedly different from the benchmark at all the locations  It might be attributed to the fact that in the winter, the relatively colder temperature does not affect all the zones, while a transition in temperature, as observed in April, may play a decisive role in operating conditions in almost all the zones. Zones 1 and 3 never figure in the list and no marked deviation from the benchmark in any climate condition (non-rejection of the null hypothesis at all seasons) is observed. This homogeneity might be attributed to the fact that these two zones correspond to a relatively warmer climate and hence climate dependence on the operating conditions are not present here. Although, we have to suppress the zone identity for confidentiality issue, the findings are as corroborated by the domain knowledge experts. \\

To conclude this section, we emphasize the appropriateness of the model-based bias correction approach. We try to explore a `what-if' type scenario and try to assess the validity of the findings if we assume the mileages at failure in each segment to be normally distributed, instead of the exponential assumption. Our main objective remains same i.e to identify the significantly different segments with respect to MMTF values. If we assume that the mileages at failure for the $i$-th segment be normally distributed with mean $\theta_i$ and variance $\sigma_i^2$, the testing problem is still the same as in (\ref{eq:datahyp}). We perform single sample both-sided t-test for each of the segments and obtain the array of $p$-values over all the segments. Computation of robust estimates may be done similarly as mentioned, but for the bias corrected estimators we follow algorithms given in Cheng et al. (2015) (for $\hat{\pi}_0^U$) and Biswas (2019) (for $\hat{\pi}_0^E$) in stead of Algorithms 1 and 2 for obvious reasons. The estimates of $\pi_0$ under normality assumption are reported in Table 6. 
\begin{singlespacing}
\begin{table}[H]\label{tab:datapi0estim}
\caption{Estimate of $\pi_0$ using different estimators under normal model assumption for the synthetic data.}
\center
\begin{tabular}{cccccc}
\hline
\\
$\hat{\pi}_0^B$ & 0.4961 & $\hat{\pi}_0^L$ & 0.4559 & $\hat{\pi}_0^A$ & 0.5111\\
\\
$\hat{\pi}_0^P$ & 0.5707 & $\hat{\pi}_0^H$ & 0.4497 & $\hat{\pi}_0^D$ & 0.5773\\
\\
$\hat{\pi}_0^S$ & 0.6786 & $\hat{\pi}_0^U$ & 0.0000 & $\hat{\pi}_0^E$ & 0.0000\\
\\
\hline
\end{tabular}
\end{table}
\end{singlespacing}
The robust estimators are seen to underestimate $\pi_0$, while the bias corrected estimators get disrupted owing to the inappropriate model assumption and hence misleading effect size of test, upper tail probability and expectation of non-null $p$-values. Thus, appropriate model-based bias correction seems to be appropriate and efficient by bringing out more power in adaptive algorithms, while the findings may be misleading        when not applied with adequate confidence on underlying model assumption. As a result, the necessary modification of bias correction technique under exponential model seems to be only way out, particularly while dealing with multiple testing problem arising from segmented failure data, usually encountered in survival and reliability studies.

\section{Discussion}
We have approached the problem of estimating $\pi_0$ and thus construction of adaptive FDR-controlling procedure from suitable model assumptions and a common test for all the hypotheses to be tested. Within the framework suggested in Cheng et al. (2015) and Biswas (2019), we have tried to develop methods for estimation of $\pi_0$ under exponential model and presented a simple adaptive Benjamini-Hochberg algorithm in a spirit similar to Craiu and Sun (2008), which is shown to be more efficient than its counterparts for simulated as well as real life synthetic data. The current work also motivates the Storey's bootstrap estimator for $\pi_0$ and the $\pi_0$-estimator based on sum of all $p$-values based on $P(p\geq V)$. The cases of $V$ being degenerated at some $\lambda$ and $V$ being uniformly distributed over $(0,1)$ have also been discussed. This may motivate other choices of $V$ for further study of model based $\pi_0$ estimators. The study on $V$ having negatively skewed density function over $(0,1)$ is presently under consideration, which tries to give more importance to the $p$-values corresponding to true null hypotheses and the construction of new estimators in future. Though the results presented in the current work strengthens the foundations of bias-corrected estimation of $\pi_0$ in general, the distinguishing feature of this work lies on the innovative application of multiple testing procedure to segmented failure data. To the best of our knowledge, such procedures have never been applied to answer such interesting research questions framed in section 6 related to large scale industrial data. In this work, however, we have focused on presenting and motivating a simple yet powerful technique of identifying significantly different segments in terms of the performance of automobile and exploring the effect of zone of operation coupled with climate, that too under exponential model assumption. The synthetic data explored in this work pose several other issues that may be solved by the application of modified methods, which are to be formulated in future.\\

This analysis of the real life synthetic data is based on one year data and may be carried out on the basis of monthly or even weekly data associated with the component failures. Owing to limited number of such failure data in each segment,  one has to use the standard failure models like exponential or Weibull. Instead, if one uses the usual Gaussian model to describe the failure pattern, then one is expected to commit a gross mistake and consequently, a false perception on the MMTF may be reached. This issue has been addressed with the same failure data. Instead of the exponential model, the normal probability model has been used and the test for equality of respective means in all 84 segments with the same null hypothesis point representing the bench mark, as being done in the usual multiple testing procedure, has been attempted. Interestingly, the test for normality at majority of  the 84 segments fails miserably and hence conclusion on the basis of the test for MMTF with reference  to the benchmark under normality will give a wrong signal about the true status of MMTF in those segments.\\

Usually, in all studies related to life time data analysis, the time to event (failure) is taken as the response variable. Although various  time to event distribution in the multivariate setup are available in the literature (Marshall and Olkin 2007; Crowder 2016), all the variables are measured with respect to the same scale. However, in a two-dimensional warranty data base the two variables are measured in a different scale and hence special treatments to analyze such data are called for. Moreover, even if the data on the two variables are available, standard warranty analysis fail to exploit both of them in a systematic manner, by neglecting the usage data. Wu and Meeker (2002) has considered the number of failures, progressively counted over time among the units which have been manufactured in immediately preceding time points and adopted a group sequential procedure (Jennison and Turnbull 1999) to address the emerging issue identification problem. They have only recorded the `month in service' (MIS) of the vehicles and have not used the individual mileage data. Instead, one may use the technique suggested here for progressively collected weekly or monthly data for all the components at all locations, available through warranty data base. The method suggested here is expected to perform equivalently, if not better, as the mileage at failure data is used to assess the performance of several components progressively, over time at various locations.\\

As the proposed method makes assumption regarding the distribution of the mileage to failure data, we should accept the fact that, the proposed estimator is not universally suitable in all situations. At the same time, multiple testing problem in a non-Gaussian framework seems to be novel and may cover all parametric models for scenarios where non-negative valued random variables seem to be appropriate. In such a framework, we have introduced two simple estimators for $\pi_0$ which simultaneously reduces the bias and variance of the existing estimator over a relatively important part of the parameter space. The behaviour of such estimator is studied through extensive simulation studies and the new estimator is shown to be more precise under some practical assumptions in comparison to those available in the existing literature. Involvement of numerical or Monte-Carlo integration for each segment makes the proposed method rather computation intensive. This extra labour is expected to be compensated by the gain in precision of the analysis, thus meaningfully addressing the multiple testing problem in a non-Gaussian set up. 
 \section*{Acknowledgement}
 The third author (AC) acknowledges the contribution of Mr. Soumen De, formerly of General Motors Tech Center India (GMTCI), Bengaluru, India for his association with the previous works based on the synthetic data, used here.
\section*{References}
\begin{small}
Benjamini, Y., and Hochberg, Y. (1995), ``Controlling the false discovery rate: a practical and powerful approach to multiple testing," \textit{Journal of the royal statistical society. Series B (Methodological)}, 289-300.\\
\\
Benjamini, Y., and Hochberg, Y. (2000), ``On the adaptive control of the false discovery rate in multiple testing with independent statistics," \textit{Journal of educational and Behavioral Statistics}, 25(1), 60-83.\\
\\
Biswas, A. (2019), ``Estimating Proportion of True Null Hypotheses based on Sum of p-values and application to microarrays," \textit{arXiv preprint arXiv:1904.13282}.\\
\\
Chen, J. J., Wang, S. J., Tsai, C. A., and Lin, C. J. (2007), ``Selection of differentially expressed genes in microarray data analysis," \textit{The pharmacogenomics journal}, 7(3), 212-220.\\
\\
Cheng, Y., Gao, D., and Tong, T. (2015), ``Bias and variance reduction in estimating the proportion of true-null hypotheses," \textit{Biostatistics}, 16(1), 189-204.\\
\\
Craiu, R. V., and Sun, L. (2008), ``Choosing the lesser evil: trade-off between false discovery rate and non-discovery rate," \textit{Statistica Sinica,} 18, 861-879.\\
\\
Crowder, M. J. (2012), ``Multivariate survival analysis and competing risks," \textit{Chapman and Hall/CRC}.\\
\\
Gianetto, Q. G., Combes, F., Ramus, C., and Gianetto, M. Q. G. (2019), Package ‘cp4p’.\\
\\
Guan, Z., Wu, B., and Zhao, H. (2008), ``Nonparametric estimator of false discovery rate based on bernštein polynomials," \textit{Statistica Sinica, }18, 905-923.\\
\\
Gupta, S. K., De, S., and Chatterjee, A. (2014), ``Warranty forecasting from incomplete two-dimensional warranty data," \textit{Reliability Engineering \& System Safety, }126, 1-13.\\
\\
Gupta, S. K., De, S., and Chatterjee, A. (2017), ``Some reliability issues for incomplete two-dimensional warranty claims data," \textit{Reliability Engineering \& System Safety, }157, 64-77.\\
\\
Hung, H. J., O'Neill, R. T., Bauer, P., and Kohne, K. (1997), ``The behavior of the p-value when the alternative hypothesis is true," \textit{Biometrics}, 11-22.\\
\\
Jennison, C., and Turnbull, B. W. (1999), ``\textit{Group Sequential Methods With Applications to Clinical Trials}," FL Chapman \& Hall/CRC.\\
\\
Jiang, H., and Doerge, R. W. (2008), "Estimating the proportion of true null hypotheses for multiple comparisons," \textit{Cancer informatics}, 6, 117693510800600001.\\
\\
Langaas, M., Lindqvist, B. H., and Ferkingstad, E. (2005), ``Estimating the proportion of true null hypotheses, with application to DNA microarray data,"\textit{ Journal of the Royal Statistical Society: Series B (Statistical Methodology), }67(4), 555-572.\\
\\
Lawless, J. F. (1998), ``Statistical analysis of product warranty data," \textit{International Statistical Review, }66(1), 41-60.\\
\\
Marshall, A. W., and Olkin, I. (2007), ``\textit{Life distributions}" (Vol. 13). Springer, New York.\\
\\
Nettleton, D., Hwang, J. G., Caldo, R. A., and Wise, R. P. (2006), ``Estimating the number of true null hypotheses from a histogram of p values," \textit{Journal of agricultural, biological, and environmental statistics, }11(3), 337.\\
\\
Ostrovnaya, I., and Nicolae, D. L. (2012), ``Estimating the proportion of true null hypotheses under dependence," \textit{Statistica Sinica}, 1689-1716.\\
\\
Pounds, S., and Cheng, C. (2006), ``Robust estimation of the false discovery rate," \textit{Bioinformatics, }22(16), 1979-1987.\\
\\
Storey, J. D. (2002), ``A direct approach to false discovery rates," \textit{Journal of the Royal Statistical Society: Series B (Statistical Methodology), }64(3), 479-498.\\
\\
Storey, J. D., and Tibshirani, R. (2003), ``SAM thresholding and false discovery rates for detecting differential gene expression in DNA microarrays,"  \textit{The analysis of gene expression data} (pp. 272-290). Springer, New York, NY.\\
\\
Wang, H. Q., Tuominen, L. K., and Tsai, C. J. (2010), ``SLIM: a sliding linear model for estimating the proportion of true null hypotheses in datasets with dependence structures,"\textit{ Bioinformatics}, 27(2), 225-231.\\
\\
Wu, H., and Meeker, W. Q. (2002), ``Early detection of reliability problems using information from warranty databases," \textit{Technometrics, }44(2), 120-133.\\

\end{small}

\end{document}